\magnification=\magstep1

\catcode`@=11
\def\displaylines#1{\displ@y
  \halign{\hbox to\displaywidth{$\@lign\hfil\displaystyle{}##\hfil$}\crcr
    #1\crcr}}
\catcode`@=12

\newcount\contaparag
\newcount\contaeq

\let\TeXeqno\eqno
\let\TeXeqalignno\eqalignno

\def\eqno#1{\global\advance\contaeq 1
  \ifx#1\undefined
  \xdef#1{\bgroup\noexpand\rm
    (\the\contaparag.\the\contaeq)\egroup}%
  \wlog{\def\string#1{(\the\contaparag.\the\contaeq)}}%
  \else
    \def\gianota{#1}%
    \ifx\gianota\space
    \else
      \messaggiogianota
    \fi
  \fi
  \TeXeqno#1}

\def\messaggiogianota{\wvd{}\wvd{ONE HAS ALREADY GIVEN A MEANING TO
  \expandafter\string\gianota\space }\wvd{}}

\def\wvd#1{\immediate\write 16 {#1}}

\let\TeXeqalignno\eqalignno

\def\eqalignno#1{\defs#1\finedefs{#1}}

\def\defs#1&#2&#3\cr#4{%
  \def\arg{#3}%
  \ifx#3\undefined
    \global\advance\contaeq 1
    \xdef#3{\bgroup\noexpand\rm
      (\the\contaparag.\the\contaeq)\egroup}%
    \wlog{\def\string#3{(\the\contaparag.\the\contaeq)}}%
  \else
    \def\gianota{#3}%
    \ifx\gianota\space
    \else
      \messaggiogianota
    \fi
  \fi
  \ifx#4\finedefs
    \let\next\finedefs
  \else
    \def\next{\defs#4}%
  \fi
  \next}

\let\finedefs\TeXeqalignno

\def\accorpa #1#2{%
  {%
   \rm\expandafter\expandafter\expandafter\ccrp\expandafter #1#2}}

\def\ccrp\bgroup\rm(#1.#2)\egroup
  \bgroup\rm(#3.#4)\egroup{(#1.#2-#4)}

% A ogni paragrafo nuovo usare \nuovoparagrafo seguito dal numero del
% paragrafo prima di cio' che c'e' gia' scritto (che andra' a capo
% senza righe vuote intermedie)

\def\nuovoparagrafo #1 {\contaparag #1 \contaeq 0 }

% numerazione dei riferimenti bibliografici

\newcount\contaref

\def\numeraref#1{%
  \def\argref{#1}%
  \ifx\argref\fineelenco
    \let\next\relax
  \else
    \advance\contaref 1
    \edef#1{\the\contaref}%
    \wlog{\def\string#1{#1}}%
    \let\next\numeraref
  \fi
  \next}

\def\fineelenco{\fineelenco}

%%%%%%%%%%%%%%%%%%%%%%%%%%%%%%%%%%%%%%%%%%%%%%%%%%%%%%%%%%%%%

\topskip=14pt
\mathsurround=2pt
\def\wbox#1/{\vbox{\hrule\hbox{\vrule height#1\kern#1\vrule height#1}\hrule}}
\newbox\boxfine \setbox\boxfine=\hbox{\wbox6pt/}
\def\fine{\ \kern6pt\copy\boxfine}
\def\media{\mkern12mu\hbox{\vrule height4pt depth-3.2pt width5pt} \mkern-16.5mu\int}
\mathchardef\TeXchi="011F
\font\maiuscolo=cmcsc10
\font\maiuscolino=cmcsc10 at 8 pt
\font\piccolo=cmr12 at 8 pt
\font\piccolobf=cmbx10 at 7pt
\font\piccolorm=cmr12 at 7pt

\font\tito=cmbx10 scaled \magstep2

\font\tenmsb=msbm10  
\font\macc=cmtt10
\font\sevenmsb=msbm7
\font\fivemsb=msbm5
\newfam\msbfam
\textfont\msbfam=\tenmsb
\scriptfont\msbfam=\sevenmsb
\scriptscriptfont\msbfam=\fivemsb
\def\Bbb{\fam\msbfam }
 \def\enne{{\Bbb N}}  
\def\erre{{\Bbb R}}

\def\chi{{\setbox0 =\hbox{$\mathsurround=0pt\TeXchi$}\hbox
  {\raise\dp0 \copy0 }}}

\def\teta{\vartheta}
\def\ti{\tau_i}
\def\eps{\varepsilon}

\def\mez{{1\over 2}}
\def\quarto{{1\over 4}}

\def\laplconf{L_{g}}

\def\qi{q_i}
\def\ui{u_i}
\def\vi{v_i}
\def\yi{y_i}
\def\privato{\,\backslash\,}
\def\interno#1{\mathop{#1}\limits^{\circ}}
\numeraref
\alm
\changyang
\cherrier
\dma
\escobaruno
\escobartre
\escobardue
\gilbarg
\hamza
\hanliuno
\hanlidue
\hu
\kellogg
\liuno
\lidue
\liliu
\lizhuuno
\lizhudue
\louzhu
\ma
\nirenberg
\ahmedou
\schoen
\schoenyau
\schoenzhang
\fineelenco

\nopagenumbers
\quad
\vskip2truecm\noindent
\nuovoparagrafo 1
\medskip\noindent
\centerline{\tito Compactness results in conformal  deformations}
\medskip\noindent
\centerline{\tito of Riemannian metrics  on manifolds}
\medskip\noindent
\centerline{\tito  with boundaries}
\smallskip\noindent
\vskip1truecm
\centerline{\maiuscolo Veronica Felli and Mohameden Ould Ahmedou}
\vskip1truecm\rm
\centerline{Scuola Internazionale Superiore di Studi Avanzati (S.I.S.S.A.)}
\centerline{via Beirut, 2-4}
\centerline{34014 Trieste, Italy}
\vskip1truecm\rm
\centerline{E-mail:\quad \macc felli@sissa.it;\quad ahmedou@sissa.it}
\vskip2truecm\noindent
\centerline{\vbox{\hsize=14truecm{\bf \noindent Abstract.}\quad{\piccolorm This paper is devoted to the study of a problem arising from a geometric context, namely the conformal deformation of a Riemannian metric to a scalar flat one having constant mean curvature on the boundary. By means of blow-up analysis techniques and the Positive Mass Theorem, we show that on locally conformally flat manifolds with umbilic boundary all metrics stay in a compact set with respect to the ${\scriptstyle C^2}\!$-norm and the total Leray-Schauder degree of all solutions is equal to ${\scriptstyle-1}$. Then we deduce from this compactness result the existence of at least one solution to our problem.\medskip\noindent{\piccolobf MSC classification:}\quad 35J60, 53C21, 58G30.}}}
\vskip2truecm\noindent
\centerline{\bf 1.\quad Introduction}
\bigskip
\rm Let $(M,g)$ be an $n\!$-dimensional compact smooth Riemannian manifold with boun\-da\-ry. For $n=2$, the well-known Riemann Mapping Theorem states that an open simply connected proper subset of the plane is conformally diffeomorphic to the disk. In what can be seen as a tentative of generalization of the above problem, J. Escobar [\escobaruno] asked if $(M,g)$ is conformally equivalent to a manifold that has zero scalar curvature and whose boundary has a constant mean curvature. \par
Setting $\tilde g=u^{4\over n-2}g$ conformal metric to $g$, the above problem is equivalent to find a smooth positive solution $u$ to the following nonlinear boundary value problem on $(M,g)$: 
$$\left\{\eqalign{&-\Delta_gu+{n-2\over 4(n-1)}R_gu=0,\quad u>0,\,\quad\hbox{in}\ \interno M,\cr
&{\partial u\over\partial\nu}+{n-2\over 2}h_gu=cu^{n\over n-2},\ \ \ \qquad\qquad\quad\hbox{on}\ \partial M,\cr}\right.\eqno({\cal P})$$
where $\interno M=M\privato\partial M$ denotes the interior of $M$, $R_g$ is the scalar curvature, $h_g$ is the mean curvature of $\partial M$, $\nu$ is the outer normal with respect to $g$, and $c$ is a constant whose sign is uniquely determined by the conformal structure of $M$. Solutions of equation~$({\cal P})$ correspond, up to some positive constant, to critical points of the following function~$J$ defined on $H^1(M)\privato\{0\}$
$$J(u)={\displaystyle\,\int_M\bigg(|\nabla_gu|^2+{n-2\over \strut{4(n-1)}}R_gu^2\bigg)\,dV_g+{n-2\over\strut2}\int_{\partial M}h_g u^2\,d\sigma_g\,\over \displaystyle\,\left(\int_{\partial M}|u|^{2(n-1)\over n-2}\,d\sigma_g\right)^{\strut{n-2\over n-1}}\,}.$$
The exponent ${2(n-1)\over n-2}$ is critical for the Sobolev trace embedding $H^1(M)\hookrightarrow~L^q(\partial M)$. This embedding being not compact, the functional $J$ does not satisfy the Palais-Smale condition. For this reason standard variational methods cannot be applied to find critical points of $J$.
\headline{\ifodd\folio\hfil\maiuscolino Compactness results in conformal  deformations on manifolds with boundaries\hfil\piccolo\folio\else\piccolo\folio\hfil\maiuscolino V. Felli and M. Ould Ahmedou\hfil \fi}
\par
The regularity of the $H^1\!$-solutions of $({\cal P})$ was established by Cherrier [\cherrier], and exi\-sten\-ce results in many cases were obtained by Escobar, see [\escobaruno, \escobardue]. Related problems regarding conformal deformations of Riemannian metrics on manifolds with boundaries have been studied in [\alm, \changyang, \dma, \escobartre, \hamza, \hanliuno, \hanlidue, \lidue, \liliu, \ma, \ahmedou]; see also the references therein. 
\par
To describe our results concerning problem $({\cal P})$, we need the following notation. We use $L_g$ to denote $\Delta_g-{(n-2)/[4(n-1)]}R_g$, $B_g$ to denote ${\partial/\partial \nu}+{(n-2)/ 2}\,h_g$. Let $H$ denote the second fundamental form of $\partial M$ in $(M,g)$ with respect to the inner normal; we denote its traceless part part by $U$:
$$U(X,Y)=H(X,Y)-h_gg(X,Y).$$
\medskip\noindent
{\bf Definition 1.1.}\quad \it A point $p\in\partial M$ is called an umbilic point if $U=0$ at $p$. The boundary of $M$ is called umbilic if every point of $\partial M$ is umbilic.
\medskip\noindent
{\bf Remark 1.2.}\quad \it The notion of umbilic point is conformally invariant, namely, if $p\in\partial M$ is an umbilic point with respect to $g$, it is also an umbilic point with respect to the metric $\tilde g=\psi^{4\over n-2}g$, for any positive smooth function $\psi$ on $M$. 
\medskip\noindent\rm
Let $\lambda_1(L)$ denote the first eigenvalue of
$$\left\{\eqalign{&-L_g\varphi =\lambda\varphi,\quad\hbox{in}\ \interno M,\cr
&B_g\varphi=0,\quad\; \; \ \quad\hbox{on}\ \partial M,\cr}\right.\eqno(E_1)$$
and $\lambda_1(B)$ denote the first eigenvalue of the problem
$$\left\{\eqalign{&L_gu=0,\;\, \; \quad\hbox{in}\ \interno M,\cr
                  &B_gu=\lambda u,\quad\hbox{on}\ \partial M.\cr}\right.\eqno(E_2)$$
It is well-known (see [\escobaruno]) that the signs of $\lambda_1(B)$ and $\lambda_1(L)$ are the same and they are conformal invariants. 
\medskip\noindent
{\bf Definition 1.3.}\quad \it We say that a manifold is of positive (respectively negative, zero) type if $\lambda_1(L)>0$ (respectively $<0$, $=0$).
\medskip\noindent\rm
In this paper, we give some existence and compactness results concerning $({\cal P})$. We first describe our results for manifolds of positive type.\par
Let $(M,g)$ be a manifold of positive type. We consider the following problem
$$\left\{\eqalign{&-L_gu=0,\quad u>0,\quad\hbox{in}\ \interno M,\cr
&B_gu=(n-2)u^{n\over n-2},\ \ \  \;\hbox{on}\ \partial M.\cr}\right.\eqno({\cal P}_+)$$
Let ${\cal M}^+$ denote the set of solutions of $({\cal P}_+)$. Then we have
\medskip\noindent
{\bf Theorem 1.4.}\quad\it For $n\geq3$, let $(M,g)$ be a smooth compact $n\!$-dimensional locally conformally flat Riemannian manifold of positive type with umbilic boundary. Then  ${\cal M}^+\not=\emptyset$. Furthermore, if $(M,g)$ is not conformally equivalent to the standard ball, then there exists $C=C(M,g)$ such that for all $u\in{\cal M}^+$ we have
$${1\over C}\leq u(x)\leq C,\quad\forall\, x\in M; \quad \hbox{and}\quad\|u\|_{C^2(M)}\leq C,$$
and the total Leray-Schauder degree of all solutions to $({\cal P}_+)$ is $-1$.
\medskip\rm
Let us remark that the existence of solutions to $({\cal P}_+)$ under the condition of Theo\-rem~1.4 was already established by Escobar in [\escobaruno], among other existence results. He obtained, using the Positive Mass Theorem of Schoen-Yau [\schoenyau], that the infimum of $J$ is achieved. See also [\ahmedou] for the existence of a solution to $({\cal P}_+)$ of higher energy and higher Morse index. What is new in Theorem 1.4 is the compactness part. In fact we establish a slightly stronger compactness result. Consider, for $1<q\leq{n\over n-2}$,
$$\left\{\eqalign{&-L_gu=0,\quad u>0,\quad\hbox{in}\ \interno M,\cr
&B_gu=(n-2)u^q,\quad\;\ \ \,\; \hbox{on}\ \partial M.\cr}\right.\eqno({\cal P}_q^+)$$
Let ${\cal M}^+_q$ denote the set of solutions of $({\cal P}_q^+)$ in $C^2(M)$. We have the following 
\medskip\noindent
{\bf Theorem 1.5.}\quad \it For $n\geq3$,  let $(M,g)$ be a smooth compact $n\!$-dimensional locally conformally flat Riemannian manifold of positive type with umbilic boundary. We assume that $(M,g)$ is not conformally equivalent to the standard ball. Then there exist $\delta_0~=~\delta_0(M,g)>0$ and $C=C(M,g)>0$ such that for all $u\in\bigcup_{1+\delta_0\leq q\leq{n\over n-2}}{\cal M}^+_q$ we have
$${1\over C}\leq u(x)\leq C,\quad \forall\, x\in M,\quad\hbox{and}\quad \|u\|_{C^2(M)}\leq C.$$
\medskip\noindent
\rm To prove Theorems 1.4 and 1.5 we establish compactness results for all solutions of $({\cal P}_q^+)$ and then show that the total degree of all solutions to $({\cal P}_+)$ is $-1$. To do this we perform some fine blow-up analysis of possible behaviour of blowing-up solutions of $({\cal P}_q^+)$ which, together with the Positive Mass Theorem by Schoen and Yau~[\schoenyau] (see also [\escobartre]), implies energy independent estimates for all solutions of $({\cal P}_q^+)$. \par
When $(M,g)$ is a $n\!$-dimensional ($n\geq3$) locally conformally flat manifold without boundary, such compactness results based on blow-up analysis and energy independent estimates were obtained by Schoen [\schoen] for solutions of
$$-L_gu=n(n-2)u^q,\quad u>0,\quad \hbox{in}\ M,$$
where $1+\eps_0<q<{n+2\over n-2}$. In the same paper [\schoen] he also announced, with indications on the proof, the same results for general manifolds. Along the same approach initiated by Schoen, Z. C. Han and Y. Y. Li [\hanliuno] obtained similar compactness and existence results for the so-called Yamabe like problem on compact locally conformally flat manifolds with umbilic boundary. Other compactness results on Yamabe type equations on three dimensional Riemannian manifolds were obtained by Y. Y. Li and M. J. Zhu [\lizhudue]. \par
Now we present similar existence and compactness results for manifolds of negative type. Let $(M,g)$ be a compact $n\!$-dimensional Riemannian manifold of negative type. Consider for $1<q\leq {n\over n-2}$
$$\left\{\eqalign{&-L_gu=0,\quad u>0,\quad\hbox{in}\ \interno M,\cr
&B_gu=-(n-2)u^q,\quad\;\ \hbox{on}\ \partial M.\cr}\right.\eqno({\cal P}_q^-)$$
Let ${\cal M}^-_q$ denote the set of solutions of $({\cal P}_q^-)$ in $C^2(M)$ and ${\cal M}^-={\cal M}^-_{n\over n-2}$. We have the following
\medskip\noindent
{\bf Theorem 1.6.}\quad\it For $n\geq3$,  let $(M,g)$ be a smooth compact $n\!$-dimensional Riemannian manifold of negative type with boundary. Then ${\cal M}^-\not=\emptyset$. Furthermore, there exist $\delta_0~=~\delta_0(M,g)$ and $C=C(M,g)>0$ such that for all $u\in\bigcup_{1+\delta_0\leq q\leq{n\over n-2}}{\cal M}_q^-$
$${1\over C}\leq u(x)\leq C,\quad\forall\, x\in M;\quad \|u\|_{C^2(M)}\leq C,$$
and the total degree of all solutions of $({\cal P}_q^-)$ is $-1$.
\medskip\noindent\rm
Let us notice that apriori estimates in the above Theorem are due basically to some non\-exi\-sten\-ce Liouville-type Theorems for the limiting equations.\par
The remainder of the paper is organized as follows. In section 2 we provide the local blow-up analysis. In section 3 we establish the compactness part in Theorems 1.4 and 1.5. In section 4 we prove existence part of Theorem 1.4 while section 5 is devoted to the proof of Theorem 1.6. Finally, we collect some technical lemmas and well-known results in the appendix.

\bigskip\noindent
\centerline{\bf Acknowledgements}
\medskip
The authors are grateful to Prof. A. Ambrosetti for his interest in their work and for his generous support. V. F. is supported by M.U.R.S.T. under the national project ``Variational Methods and Nonlinear Differential Equations'' and M. O. A. research is supported by a S.I.S.S.A. postdoctoral fellowship.

\vskip2truecm
\nuovoparagrafo 2

\centerline{\bf 2.\quad Local blow-up analysis}
\bigskip\noindent
In the following, we give the definitions of isolated and isolated simple blow-up, which were first introduced by R. Schoen, see [\schoen], and adapted to the framework of boundary value problems by Y. Y. Li [\lidue], see also [\hanliuno]. 
\medskip\noindent 
{\bf Definition 2.1.}\quad \it Let $(M,g)$ be a smooth compact $n\!$-dimensional Riemannian manifold with boundary, and let $\bar r>0$, $\bar c>0$, $\bar x\in\partial M$, $f\in C^0(\overline{B_{\bar r}(\bar x)})$ be some positive function where  $B_{\bar r}(\bar x)$ denotes the geodesic ball in $(M,g)$ of radius $\bar r$ centered at $\bar x$. Suppose that, for some sequences $q_i={n\over n-2}-\tau_i$, $\tau_i\to 0$, $f_i\to f$ in $C^0(\overline{B_{\bar r}(\bar x)})$, $\{u_i\}_{i\in \enne}$ solves
$$\left\{\eqalign{&-\laplconf u_i=0,\quad\ui>0,\quad \hbox{in}\ B_{\bar r}(\bar x),\cr
                  &B_gu_i=(n-2)f_i^{\tau_i}\ui^{\qi},\quad\hbox{on}\ \partial M\cap B_{\bar r}(\bar x).\cr}\right.\eqno\dueuno_i$$ 
We say that ${\bar x}$ is an isolated blow-up point of $\{\ui\}_i$ if there exists a sequence of local maximum points $x_i$ of $\ui$ such that $x_i\to \bar x$ and, for some $C_1>0$,
$$\lim_{i\to\infty}\ui(x_i)=+\infty\quad\hbox{and}\quad\ui(x)\leq  C_1 d(x,x_i)^{-{1\over q_i-1}},\quad\forall\, x\in B_{\bar r}(x_i),\ \forall\, i.$$
\medskip\rm
To describe the behaviour of blowing-up solutions near an isolated blow-up point, we define spherical averages of $\ui$ centered at $x_i$ as follows
$$\bar \ui(r)=\media_{M\cap\partial B_r(\bar x)}\ui={1\over {\rm Vol}_g(M\cap\partial B_r(\bar x))}\int_{M\cap\partial B_r(\bar x)}\ui.$$
Now we define the notion of isolated blow-up point.
\medskip\noindent
{\bf Definition 2.2.}\quad \it Let $x_i\to \bar x$ be an isolated blow-up point of $\{\ui\}_i$ as in Definition 2.1. We say that $x_i\to \bar x$ is an isolated simple blow-up point of $\{\ui\}$ if, for some positive constants $\tilde r\in(0,\bar r)$ and $C_2>1$, the function $\bar w_i(r):=r^{1\over q_i-1}\bar\ui(r)$ satisfies, for large $i$,
$$\bar w_i'(r)<0\quad\hbox{for}\ r\ \hbox{satisfying}\quad  C_2\ui^{1-q_i}(x_i)\leq r\leq\tilde r.$$
\medskip\noindent\rm
Let us introduce the following notation
$$\eqalign{&\erre^n_+=\{(x',x^n)\in \erre^{n-1}\times \erre:\ x^n>0\},\quad B^+_\sigma(\bar x)=\{x=(x',x^n)\in\erre^n_+:\ |x-\bar x|<\sigma\},\cr
&B^+_\sigma=B_\sigma^+(0),\quad \Gamma_1(B_\sigma^+(\bar x))=\partial B_\sigma^+(\bar x)\cap\partial\erre^n_+,\quad \Gamma_2(B_\sigma(\bar x))=\partial B_\sigma(\bar x)\cap\erre^n_+.\cr}$$
Let $\{f_i\}\subset C^1(\Gamma_1(B_3^+))$ be a sequence of functions satisfying, for some positive constant $C_3$,
$$f_i\ \mathop{\longrightarrow}\limits_{i\to\infty}\ f\quad\hbox{in}\ C^1(\Gamma_1(B_3^+)),\quad\|f_i\|_{L^\infty(\Gamma_1(B_3^+))}\leq C_3\eqno\duedue$$
where $f\in C^1(\Gamma_1(B_3^+))$ is some positive function. Suppose that $\{\vi\}_i\subset C^2(\overline{B_3^+})$ is a sequence of solutions to
$$\left\{\eqalign{& -\Delta \vi=0,\quad v_i>0,\quad\ \ \ \ \hbox{in}\ B_3^+,\cr
&{\partial \vi\over\partial x^n}=-(n-2)f_i^{\ti} \vi^{\qi},\quad\hbox{on}\ \Gamma_1(B_3^+).\cr}\right.\eqno\duetre_i$$
The following Lemma gives a Harnack inequality.
\medskip\noindent
{\bf Lemma 2.3.}\quad\it Assume \duedue\ and let $\{\vi\}_i$ satisfy $\duetre_i$. Let $0<\bar r<{1\over8}$, $\bar x\in \Gamma_1(\overline{B^+_{1/8}})$ and suppose that $x_i\to x$ is an isolated blow-up point of $\{\vi\}_i$. Then, for all $0<r<\bar r$,
$$\sup _{B^+_{2r}(x_i)\,\backslash\,B^+_{r\over2}(x_i)}\vi\leq C_4\inf_{B^+_{2r}(x_i)\,\backslash\,B^+_{r\over2}(x_i)}\vi,$$
where $C_4>0$ is some positive constant independent of $i$ and $r$.\rm
\medskip\noindent
{\bf Proof.}\quad Without loss of generality, we assume that $x_i\in\Gamma_1(\overline{B^+_{1/8}})$. For $0<r<\bar r$, let us consider 
$$\tilde \vi(y):=r^{1\over \qi-1}\vi(ry+x_i).$$
Then $\tilde\vi$ satisfies 
$$\left\{\eqalign{&-\Delta\tilde\vi=0,\quad\tilde \vi>0,\qquad\qquad\qquad\qquad\,\hbox{in}\ A_i,\cr
&{\partial\tilde\vi\over\partial y^n}=-(n-2)f_i^{\ti}(ry+x_i)\tilde\vi^{q_i-1}\tilde\vi,\quad\hbox{on}\ \Gamma_1(A_i),\cr}\right.$$
where $A_i=\left\{y\in\erre^n:\ {1\over3}<|y|<3,\ ry+x_i\in\erre^n_+\right\}$. From Definition 2.1 we know that
$$\tilde \vi\leq C_1\quad\hbox{in}\ A_i,$$
where $C_1$ depends neither on $r$ nor on $i$. In view of \duedue, from Lemma 6.1 (standard Harnack) in the appendix we obtain that for some constant $c>0$
$$\max_{\widetilde A_i}\tilde\vi\leq c\min _{\widetilde A_i}\tilde\vi$$
where $\widetilde A_i=\left\{y\in\erre^n:\ \mez<|y|<2,\ ry+x_i\in\erre^n_+\right\}$; the proof of the Lemma is thereby completed.\hfill\fine
\medskip\noindent
{\bf Lemma 2.4.}\quad \it Suppose that $\{\vi\}_i$ satisfies $\duetre_i$ and $\{x_i\}_i\subset\Gamma_1(B_1^+)$ is a sequence of local maximum points of $\{\vi\}_i$ in $\overline{B_3^+}$ satisfying
$$\{\vi(x_i)\}\quad\hbox{is bounded},$$
and, for some constant $C_5$,
$$|x-x_i|^{1\over {q_i}-1}\vi(x)\leq C_5,\quad \forall\, x\in B_3^+.\eqno\duequattro$$
Then
$$\limsup_{i\to\infty}\max_{\overline{B^+_{1/4}(x_i)}}\vi<\infty.\eqno\duecinque$$
\medskip\noindent
\rm 
{\bf Proof.}\quad By contradiction, suppose that, under the assumptions of the Lemma, \duecinque\ fails, namely that, along a subsequence, for some $\tilde x_i\in\overline{B^+_{1/4}(x_i)}$ we have
$$\vi(\tilde x_i)=\max_{\overline{B^+_{1/4}(x_i)}}\vi\ \mathop{\longrightarrow}\limits_{i\to\infty}\ +\infty.$$
It follows from \duequattro\ that $|\tilde x_i-x_i|\to 0$. Let us now consider
$$\xi_i(z)=\vi^{-1}(\tilde x_i)\vi(\tilde x_i+\vi^{1-\qi}(\tilde x_i)z)$$
defined on the set
$$B^{-T_i}_{{1\over 8}\vi^{\qi-1}(\tilde x_i)}:=\left\{z\in\erre^n:\ |z|<{1\over8}\vi^{\qi-1}(\tilde x_i)\quad\hbox{and}\quad z^n>-T_i\right\}$$
where $T_i=\tilde x_i^n\vi^{\qi-1}(\tilde x_i)$. In view of $\duetre_i$, $\xi_i$ satisfies
$$\left\{\eqalign{&-\Delta\xi_i=0,\quad\xi_i>0,\qquad z\in B^{-T_i}_{{1\over 8}\vi^{\qi-1}(\tilde x_i)},\cr
&{\partial\xi_i\over\partial z^n}=-(n-2)f_i^{\ti}\xi_i^{\qi},\quad z\in\partial B^{-T_i}_{{1\over 8}\vi^{\qi-1}(\tilde x_i)}\cap\{z=(z',z^n)\in\erre^n:\ z^n=-T_i\},\cr}\right.$$
and
$$\xi_i(z)\leq\xi_i(0)=1,\quad\forall\, z\in B^{-T_i}_{{1\over 8}\vi^{\qi-1}(\tilde x_i)}.$$
It follows from \duequattro\ that
$$|z|^{1\over \qi-1}\xi_i(z)\leq C_1,\quad\forall\, z\in B^{-T_i}_{{1\over 8}\vi^{\qi-1}(\tilde x_i)}.$$
Since $\{\xi_i\}_i$ is locally bounded, applying $L^p\!$-estimates, Schauder estimates, the Harnack inequality, and Lemma 6.1, we have that, up to a subsequence, there exists some positive function $\xi$ such that
$$\lim_{i\to\infty}\|\xi_i-\xi\|_{C^2(\erre^n_{-T_i}\cap\overline{B_R})}=0,\quad\forall\, R>1,$$
where $\erre^n_{-T_i}=\{z=(z',z^n)\in\erre^n:\ z^n>-T_i\}$ and, for $T=\lim_{i\to\infty}T_i\in[0,+\infty]$, $\xi$ satisfies
$$\left\{\eqalign{&-\Delta \xi=0,\quad\xi>0,\quad\quad\hbox{in}\ \erre^n_{-T},\cr
&{\partial\xi\over\partial z^n}=-(n-2)\xi^{n\over n-2},\quad\hbox{on}\ \partial \erre^n_{-T}.\cr}\right.\eqno\duecinqueprimo$$
Let us prove that $T<\infty$. Indeed, if we assume by contradiction that $T=+\infty$, we have that $\xi$ is a harmonic bounded function in $\erre^n$. The Liouville Theorem yields that $\xi$ is a constant and this is in contradiction with \duequattro.\par
Therefore $T<\infty$. Let us prove that $T=0$. Since problem \duecinqueprimo, up to a translation, satisfies the assumptions of the uniqueness Theorem by Li and Zhu [\lizhuuno], we deduce that $\xi$ is of the form
$$\xi(x',x^n)=\left[{\lambda\over(1+\lambda(x^n-T))^2+\lambda^2|x'-x_0'|^2}\right]^{n-2\over 2}$$
for some $\lambda>0$, $x'_0\in\erre^{n-1}$. Since $0$ is a local maximum point for $\xi$, it follows that $x_0'=0$ and $T=0$. Furthermore the fact that $\xi(0)=1$ yields $\lambda=1$. It follows that, for all $R>1$
$$\min_{\overline {B}^{-T_i}_{R\vi^{-(\qi-1)}(\tilde x_i)}(\tilde x_i)}\vi=\vi(\tilde x_i)\min_{\overline {B}_R^{-T_i}(0)}\xi_i\ \mathop{\longrightarrow}\limits_{i\to\infty}\ \infty.$$
Since $\{\vi(x_i)\}_i$ is bounded, we have that, for any $R>1$, $x_i\not\in \bar B^{-T_i}_{R\vi^{\qi-1}(\tilde x_i)}(\tilde x_i)$ for large $i$, namely 
$$R<\vi^{\qi-1}(\tilde x_i)|\tilde x_i-x_i|.$$
Hence we have that
$$|\tilde x_i-x_i|^{1\over\qi-1}\vi(\tilde x_i)>R^{1\over \qi-1}$$
which contradicts \duequattro.\hfill\fine
\medskip\noindent
{\bf Proposition 2.5.}\quad\it Let $(M,g)$ be a smooth compact $n\!$-dimensional locally conformally flat Riemannian manifold with umbilic boundary, and let $x_i\to\bar x$ be an isolated simple blow-up point of $\{\ui\}_i$. Then for any sequences of positive numbers $R_i\to\infty$, $\eps_i\to0$ there exists a subsequence $\{u_{j_i}\}_i$ (still denoted as $\{\ui\}_i$) such that
$$r_i:=R_i\ui^{1-\qi}(x_i)\ \mathop{\longrightarrow}\limits_{i\to\infty}\ 0,\quad x_i\in\partial M,$$
and
$$\eqalign{&\left\|\ui^{-1}(x_i)\ui(\exp_{x_i}(y\ui^{1-\qi}(x_i))-\left({1\over (1+y^n)^2+|y'|^2}\right)^{n-2\over2}\right\|_{C^2(B_{3R_i}(0))}\cr
&\qquad+\left\|\ui^{-1}(x_i)\ui(\exp_{x_i}(y\ui^{1-\qi}(x_i))-\left({1\over (1+y^n)^2+|y'|^2}\right)^{n-2\over2}\right\|_{H^1(B_{3R_i}(0))}<\eps_i.\cr}$$
Moreover, for all $2r_i\leq d(x,x_i)\leq\tilde r/2$,
$$\ui(x)\leq C_6\ui^{-1}(x_i)d(x,x_i)^{2-n},$$
where $C_6$ is some positive constant independent of $i$, and
$$\ui(x_i)u_i\ \mathop{\longrightarrow}\limits_{i\to\infty}\ aG(\cdot, \bar x)+b\quad\hbox{in}\ C^2_{\rm loc}(\overline{B_{\tilde r}(\bar x)}\,\backslash\,\{\bar x\})$$
where $a>0$, $b$ is some nonnegative function satisfying
$$\left\{\eqalign{&L_gb=0,\quad\hbox{in}\ B_{\tilde r}(\bar x)\,\backslash\,\{\bar x\},\cr
                  &B_gb=0,\quad\hbox{on}\ B_{\tilde r}(\bar x)\cap\partial M,\cr}\right.$$
and $G(\cdot, \bar x)$ is the Green's function satisfying
$$\left\{\eqalign{&-L_gG(\cdot,\bar x)=0,\quad\hbox{in}\ M\,\backslash\,\{\bar x\},\cr
                  &B_gG(\cdot,\bar x)=0,\quad\quad\hbox{on}\ \partial M\,\backslash\,\{\bar x\}.\cr}\right.\eqno\green$$
\rm\medskip\noindent
To prove Proposition 2.5 we need some preliminary results. Hence forward we use~$c$, $c_1,\ c_2, \dots$ to denote positive constants which may vary from formula to formula and which may depend only on $M$, $g$, $n$, and $\bar r$.
\medskip\noindent
{\bf Lemma 2.6.}\quad\it Let $x_i\to0$ be an isolated blow-up point of $\{\vi\}_i$ with $\vi$ solutions of $\duetre_i$. Then, for any $R_i\to\infty$ and $\eps_i\to0$, there exists a subsequence of $\{\vi\}_i$, still denoted by $\{\vi\}_i$, such that 
$$r_i:=R_i\vi^{1-\qi}(x_i)\ \mathop{\longrightarrow}\limits_{i\to\infty}\ 0$$
and
$$\eqalign{&\left\|\vi^{-1}(x_i)\vi(x\vi^{1-\qi}(x_i)+x_i)-\left({1\over (1+x^n)^2+|x'|^2}\right)^{n-2\over2}\right\|_{C^2(B^+_{3R_i})}\cr
&\qquad+\left\|\vi^{-1}(x_i)\vi(x\vi^{1-\qi}(x_i)+x_i)-\left({1\over (1+x^n)^2+|x'|^2}\right)^{n-2\over2}\right\|_{H^1(B^+_{3R_i})}<\eps_i.\cr}$$
\rm\medskip\noindent
{\bf Proof.}\quad Let us set
$$\tilde \vi(z)=\vi^{-1}(x_i)\vi(\vi^{1-\qi}(x_i)z+x_i), \quad z\in B^{-T_i}_{\vi^{\qi-1}(x_i)},$$
where $T_i=x_i^n\vi^{\qi-1}(x_i)$. It is clear that $\tilde\vi$ satisfies
$$\left\{\eqalign{&-\Delta \tilde\vi=0,\qquad\qquad\quad\qquad\quad\quad\qquad\qquad\hbox{in}\ B^{-T_i}_{\vi^{\qi-1}(x_i)},\cr
&{\partial \tilde\vi\over\partial z^n}=-(n-2)f_i^{\ti}(\vi^{1-\qi}(x_i)z+x_i)\tilde\vi^{\qi},\quad\hbox{on}\ \partial B^{-T_i}_{\vi^{\qi-1}(x_i)}\cap\{z\in\erre^n:\ z^n=-T_i\}.\cr}\right.$$
Let us prove that $\tilde v_i$ is uniformly bounded. By definition of isolated blow-up point, we have that
$$|z|^{1\over\qi-1}\tilde \vi(z)\leq C_1,\quad\forall\, z\in B^{-T_i}_{\vi^{\qi-1}(x_i)}.\eqno\duesei$$
It follows from \duesei, Lemma 2.3, and the Harnack inequality that $\tilde \vi$ is uniformly bounded in $B^{-T_i}_{\vi^{\qi-1}(x_i)}\cap\bar B_{R}$ for any $R>0$. Then, up to a subsequence, setting $T=\lim_{i\to\infty}T_i\in[0,+\infty]$, $\tilde\vi$ converges to some $\tilde v$ in $C^2_{\rm loc}(\erre^n_{-T})$ satisfying
$$\left\{\eqalign{&-\Delta \tilde v=0,\quad \tilde v>0,\,\quad\ \quad\hbox{in}\ \erre^n_{-T},\cr
&{\partial \tilde v\over\partial x^n}=-(n-2)\tilde v^{n\over n-2},\;\, \quad\hbox{on}\ \partial \erre^n_{-T}\quad (\hbox{if}\,\ T<\infty).\cr}\right.\eqno\duesette$$
We claim that $T<\infty$. Indeed, if we assume by contradiction that $T=+\infty$, we have that $\tilde v$ is a harmonic bounded function in $\erre^n$. By the Liouville Theorem, this implies that $\tilde v$ is a constant and this is in contradiction with \duesei.\par
Therefore $T<\infty$ and it follows from Li and Zhu uniqueness result [\lizhuuno] that $T=0$, hence
$$\tilde v(x)=\left({1\over (1+x^n)^2+|x'|^2}\right)^{n-2\over 2}.$$
So, Lemma 2.6 follows.\hfill\fine
\medskip\noindent
{\bf Lemma 2.7.}\quad \it Let $x_i\to 0$ be an isolated simple blow-up point of $\{\vi\}_i$, where $\vi$ are solutions of $\duetre_i$, and 
$$|x-x_i|^{1\over\qi-1}\vi(x)\leq C_7,\quad\forall\, x\in B_2^+,$$
for some positive constant $C_7$ and
$$\bar w_i'(r)<0,\quad\forall\, r_i\leq r\leq 2.$$
Then, for each sequence $R_i\to\infty$, there exists $\delta_i>0$, $\delta_i=O(R_i^{-1+o(1)})$ such that 
$$\vi(x)\leq C_8  \vi^{-\lambda_i}(x_i)|x-x_i|^{2-n+\delta_i},\quad \forall\, r_i\leq|x-x_i|\leq 1,$$
where $r_i=R_i\vi^{1-\qi}(x_i)$, $\lambda_i=(n-2-\delta_i)(\qi-1)-1$, and $C_8$ is some positive constant independent of $i$.
\medskip\noindent
{\bf Proof.}\quad\rm For any $x\in\{x\in\erre^n:\ r_i<|x-x_i|<2\}$, using the Harnack inequality we have that
$$|x-x_i|^{1\over\qi-1}\vi(x)\leq c \bar \vi(|x-x_i|)|x-x_i|^{1\over \qi-1}.$$
Since the blow-up is isolated simple, we have that the function at the right hand side is decreasing so that we deduce
$$|x-x_i|^{1\over\qi-1}\vi(x)\leq c\bar\vi(r_i)r_i^{1\over\qi-1}$$
for some positive constant $c$. Since 
$$\bar\vi(r_i)={1\over |\Gamma_2(B^+_{r_i})|}\int_{\Gamma_2(B^+_{r_i})}\vi$$
from Lemma 2.6 we deduce that for any $r_i<|x-x_i|<2$
$$|x-x_i|^{1\over\qi-1}\vi(x)\leq R_i^{{2-n\over 2}+o(1)}$$
which yields
$$\vi^{\qi-1}(x)\leq c|x-x_i|^{-1}R_i^{{2-n\over 2}(\qi-1)+o(1)}=c|x-x_i|^{-1}R_i^{-1+o(1)}.\eqno\dueotto$$
Set $T_i=x_i^n\vi^{\qi-1}(x_i)$. From the proof of Lemma 2.6 we know that $\lim_i T_i=0$. It is not restrictive to suppose that $x_i=(0,0,\dots, 0, x_i^n)$. Thus we have that
$$|x_i^n|=o(\vi^{1-\qi}(x_i))=o(r_i).$$
So 
$$B_1^+(0)\,\backslash\, B^+_{2r_i}(0)\subset\left\{x\in \erre^n:\ {3\over2}r_i\leq|x-x_i|\leq{3\over2}\right\}.$$
Let us apply the Maximum Principle stated in the appendix (Theorem 6.2) with
$$\eqalign{&\Omega=D_i:=B_1^+(0)\,\backslash\, B^+_{2r_i}(0),\cr
            &\Sigma=\Gamma_1(D_i)=\partial D_i\cap\partial\erre^n_+,\quad\Gamma=\Gamma_2(D_i)=\partial D_i\cap\erre^n_+,\cr
&V\equiv 0,\ \quad\qquad\quad\quad \qquad\qquad h=(n-2)f_i^{\tau_i}\vi^{\qi-1},\cr
&\psi=\vi,\ \quad\qquad\quad\quad \qquad\qquad v=\varphi_i,\cr}$$
where
$$\varphi_i(x)=M_i(|x|^{-\delta_i}-\eps_i|x|^{-\delta_i-1}x^n)+A\vi^{-\lambda_i}(x_i)(|x|^{2-n+\delta_i}-\eps_i|x|^{1-n+\delta_i}x^n)-\mez\vi(x)$$
with $M_i$, $A$, $\eps_i$, $\delta_i=O(R_i^{-1+o(1)})$ to be suitably chosen and $\lambda_i=(n-2-\delta_i)(\qi-1)-1$.\par
A straightforward calculation gives
$$\eqalign{\Delta\varphi_i(x)&=M_i|x|^{-\delta_i-2}[-\delta_i(n-2-\delta_i)+O(\eps_i)]\cr
   &\quad +|x|^{-(n-\delta_i)}A\vi^{-\lambda_i}(x_i)[-\delta_i(n-2-\delta_i)+O(\eps_i)],\quad x\in D_i\cr}$$
and, taking into account \dueotto, we have 
$$\eqalign{B\varphi_i=&M_i|x|^{-\delta_i-1}\left(\eps_i-O\left(R_i^{-1+o(1)}\right)\right)\cr
           &+A\vi^{-\lambda_i}(x_i)|x|^{-n+1+\delta_i}\left(\eps_i-O\left(R_i^{-1+o(1)}\right)\right)
,\quad\hbox{on}\ \Gamma_1(D_i).\cr}$$
Apparently we can find $0<\delta_i=O\left(R_i^{-1+o(1)}\right)$ and $0<\eps_i=O\left(R_i^{-1+o(1)}\right)$, so that
$$\Delta\varphi_i\leq0,\quad\hbox{in}\ D_i,\quad{\partial\varphi_i\over\partial x^n}+(n-2)f_i^{\tau_i}\vi^{\qi-1}\leq	0,\quad\hbox{on}\ \Gamma_1(D_i).$$
Now we check $\varphi_i\geq 0$ on $\Gamma_2(D_i)$. We have that $\Gamma_2(D_i)=\Gamma_{r_i}\cup\Gamma_2(B_1^+)$ where $\Gamma_{r_i}~=~\{x\in \erre^n_+:\ |x|=r_i\}$, $\Gamma_2(B_1^+)=\{x\in\erre^n_+:\ |x|=1\}$. On $\Gamma_{r_i}$ we have that
$$\vi(x)\leq c\vi(x_i)R_i^{2-n}\eqno\duenove$$
for some positive $c$. Choose $A$ such that
$$A\vi(x_i)R_i^{2-n+\delta_i}-c\vi(x_i)R_i^{2-n}\geq0.$$
Then by \duenove\ and for $\eps_i$ small enough we have that $\varphi_i\geq0$ on $\Gamma_{r_i}$ and taking $M_i~=~\max_{\Gamma_2(B_1^+)}\vi$ we obtain $\varphi_i\geq 0$ on $\Gamma_2(B_1^+)$. Then from Theorem 6.2 we derive that $\varphi_i\geq0$, and hence
$$\eqalignno{\vi(x)\leq &M_i(|x|^{-\delta_i}-\eps_i|x|^{-\delta_i-1}x^n)&\cr
&+A\vi^{-\lambda_i}(x_i)(|x|^{2-n+\delta_i}-\eps_i|x|^{1-n+\delta_i}x^n)\quad\forall\, x\in D_i.
&\duedieci\cr}$$
By the Harnack inequality and by the assumption that the blow-up point is isolated simple, we derive
$$M_i\leq c\bar \vi(1)\leq c\teta^{1\over\qi-1}\bar\vi(\teta)\quad\forall\,\teta\in(r_i,1).\eqno\dueundici$$
From \duedieci\ and \dueundici\ we have that
$$M_i\leq c\left\{\teta^{1\over\qi-1}\left[M_i\teta^{-\delta_i}+A\vi^{-\lambda_i}(x_i)\teta^{2-n+\delta_i}\right]\right\}$$
which implies
$$M_i\teta^{n-2-\delta_i-{1\over\qi-1}}\left(1-c\teta^{{1\over\qi-1}-\delta_i}\right)\leq cA\vi^{-\lambda_i}(x_i).$$
Choosing $\teta$ such that $1-c\teta^{2\over n-2}>1/10$, we obtain that
$$M_i\leq c\vi^{-\lambda_i}(x_i)\eqno\duedodici$$
for some constant $c>0$. The conclusion of the Lemma follows from \duedieci\ and \duedodici.\hfill\fine
\medskip
The following Lemma is a consequence of the Pohozaev identity in the appendix (see Theorem 6.3), Lemma 2.6, Lemma 2.7, and standard elliptic arguments..
\medskip\noindent
{\bf Lemma 2.8.}\quad\it $\tau_i=O\left(\vi^{-2}(x_i)\right)$. In particular $\lim_i \vi^{\ti}(x_i)=1$.
\medskip\noindent
{\bf Lemma 2.9.}\quad \it Under the same assumptions of Lemma 2.7, we have that for some positive constant $C_9>0$
$$\vi(x_i)\vi(x)\leq C_9|x-x_i|^{2-n},\quad\forall x\in B_3^+,\eqno\stima$$
and 
$$\vi(x_i)\vi\ \mathop{\longrightarrow}\limits_{i\to\infty}\ a|x|^{2-n}+b\quad\hbox{in}\ C^2_{\rm loc}(B_1^+\,\backslash\, \{0\})$$
where $a$ is a positive constant and $b\geq0$ satisfies
$$\left\{\eqalign{&-\Delta b=0,\quad\hbox{in}\ B_1^+,\cr
&{\partial b\over\partial \nu}=0,\, \quad\quad\hbox{on}\ \Gamma_1(B_1^+).\cr}\right.$$
\medskip\noindent
{\bf Proof.}\quad \rm The inequality in Lemma 2.9 for $|x-x_i|<r_i$ follows immediately from Lemma 2.6 and Lemma 2.8. Let $e\in\erre^n$, $e\in \Gamma_2(B_1^+)$, and set
$$\tilde\vi(x)=\vi^{-1}(x_i+e)\vi(x).$$
Then $\tilde\vi$ satisfies
$$\left\{\eqalign{&-\Delta\tilde\vi=0,\quad\tilde\vi>0,\, \ \ \ \quad\ \quad\qquad\quad\quad \hbox{in}\ B_2^+,\cr
                  &{\partial\tilde\vi\over\partial x^n}=-(n-2)f_i^{\ti}\vi^{\qi-1}(x_i+e)\tilde\vi^{\qi},\quad\hbox{on}\ \Gamma_1(B_2^+).\cr}\right.$$
Using Lemma 2.3 and some standard elliptic estimates, it follows, after taking a subsequence, that $\tilde\vi$ converges in $C^2_{\rm loc}(\overline{B_2^+\,\backslash\,\{0\}})$ to some positive function $\tilde v\in C^2_{\rm loc}(\overline{B_2^+\,\backslash\,\{0\}})$ satisfying       
$$\left\{\eqalign{&-\Delta\tilde v=0,\quad\hbox{in}\ B_2^+\,\backslash\,\{0\},\cr
                  &{\partial\tilde v\over\partial x^n}=0,\quad\ \ \hbox{on}\ \Gamma_1(B_2^+)\,\backslash\,\{0\},\cr}\right.\eqno\duequattordici$$
where we have used Lemma 2.7 to derive the second equation in \duequattordici. By Schwartz reflection, we obtain a function (still denoted by $\tilde v$) in $B_2$ satisfying
$$\Delta\tilde v=0,\quad\hbox{in}\ B_2\privato\{0\}.$$
So by B\"ocher's Theorem, see e.g. [\kellogg], it follows that $\tilde v(x)=a_1|x|^{2-n}+b_1$, where $a_1\geq0$, $\Delta b_1=0$, and ${\partial b_1\over\partial x^n}=0$ on $\Gamma_1(B_2^+)$. Furthermore $\tilde v$ has to be singular at $x=0$. Indeed it follows from Lemma 2.3 and some standard elliptic estimates that for $0<r<2$,
$$\lim_{i\to\infty}\vi^{-1}(x_i+e)r^{1\over\qi-1}\bar\vi(r)=r^{n-2\over2}\xi(r)$$
where $$\xi(r)=\media_{\Gamma_2(B_r^+)}\tilde v.$$ 
Therefore, it follows from the definition of isolated simple blow-up point that $r^{n-2\over2}\bar\xi(r)$ is decreasing , which is impossible if $\xi$ is regular at the origin. It follows that $a_1>0$.\par
We first establish the inequality in Lemma 2.9 for $|x-x_i|=1$. Namely, we prove that
$$\vi(x_i+e)\vi(x_i)\leq c,\eqno\duequindici$$
for some constant $c>0$. Suppose that \duequindici\ is not true, then along some subsequence we have
$$\lim_{i\to\infty}\vi(x_i+e)\vi(x_i)=\infty.$$
Multiply $\duetre_i$ by $\vi^{-1}(x_i+e)$ and integrate by parts over $B_1^+$ to obtain
$$0=\int_{B_1^+}(-\Delta\vi)\vi^{-1}(x_i+e)=\vi^{-1}(x_i+e)\int_{\partial B_1^+}{\partial\vi\over\partial\nu}.$$
Hence from the boundary condition in $\duetre_i$ we have that 
$$0=(n-2)\vi^{-1}(x_i+e)\int_{\Gamma_1(B_1^+)}f_i^{\ti}\vi^{\qi}+\vi^{-1}(x_i+e)\int_{\Gamma_2(B_1^+)}{\partial\vi\over\partial\nu}.$$
Then we have
$$\eqalignno{&\lim_{i\to\infty}\left((n-2)\vi^{-1}(x_i+e)\int_{\Gamma_1(B_1^+)}f_i^{\ti}\vi^{\qi}\right)=-\lim_{i\to\infty}\left(\int_{\Gamma_2(B_1^+)}{\partial\tilde\vi\over\partial\nu}\right)&\cr
&\quad=-\int_{\Gamma_2(B_1^+)}{\partial\tilde v\over\partial\nu}=(n-2)a_1\int_{\Gamma_2(B_1^+)}|x|^{1-n}+\int_{\Gamma_2(B_1^+)}{\partial b_1\over\partial\nu}&\cr
&\quad=(n-2)a_1|{\Gamma_2(B_1^+)}|+\int_{B_1^+}(-\Delta b_1)=(n-2)a_1|S^{n-1}_+|>0.&\duesedici\cr}$$ 
On the other hand, in view of Lemma 2.6, Lemma 2.7, and \duequindici, it is easy to check that
$$(n-2)\vi^{-1}(x_i+e)\int_{\Gamma_1(B_1^+)}f_i^{\ti}\vi^{\qi}=o(1)\vi^{-1}(x_i)\vi^{-1}(x_i+e)\ \mathop{\longrightarrow}\limits_{i\to\infty}\ 0$$
which is in contradiction with \duesedici.\par
So we have established the inequality for $|x-x_i|=1$. To establish the inequality for $r_i\leq|x-x_i|\leq3$, it is sufficient to scale the problem to reduce it to the case $|x-x_i|=1$.\par
It follows from the above that $w_i=\vi(x_i)\vi\to w$ in $C^2_{\rm loc}(B_1^+\privato\{0\})$ where $w(x)=aG(\bar x,x)+b$, for some positive constant $a$ and a function $b\geq0$ satisfying
$$\left\{\eqalign{&\Delta b=0,\quad\hbox{in}\ B_1,\cr
                  &{\partial b\over\partial \nu}=0,\quad\hbox{on}\ \Gamma_1(B_1^+).\fine\cr}\right.$$
\medskip\noindent
{\bf Proof of Proposition 2.5.}\quad Since $M$ is locally conformally flat and the boundary of $M$ is umbilic, we can find a diffeomorphism $\varphi:\ B_2^+\to B_{\bar r}(\bar x)$ and $f\in C^2(\overline{B_2^+})$ some positive function such that $\varphi(0)=\bar x$ and $\varphi^*g=f^{4\over n-2}g_0$, where $g_0$ is the flat metric in $B_2^+$. Let $\vi=f \ui\circ\varphi$. It follows from the conformal invariance of $L_g$ and $B_g$ that $\vi$ satisfies equation $\duetre_i$. So the proof of Proposition 2.5 can be easily deduced from Lemma 2.6 and Lemma 2.9.~\hfill\fine
\medskip\noindent
{\bf Proposition 2.10.}\quad\it Let $(M,g)$ be a smooth compact $n\!$-dimensional locally conformally flat Riemannian manifold with umbilic boundary and $x_i\to\bar x$ be an isolated blow-up point of $\{\ui\}_i$, where $\ui$ are solutions of $\dueuno_i$. Then it is necessarily an isolated simple blow-up point.
\medskip\rm
Due to the conformal invariance of $L_g$ and $B_g$, the proof of Proposition 2.10 is reduced to the proof of the following
\medskip\noindent
{\bf Proposition 2.11.}\quad\it Let $x_i\to 0$ be an isolated blow-up point of $\{\vi\}_i$, where $v_i$ are solutions of $\duetre_i$. Then it is an isolated simple blow-up point.
\medskip\noindent
{\bf Proof.}\quad \rm It follows from Lemma 2.6 that
$$\bar w_i'(r)<0\quad\hbox{for every}\quad C_2\vi^{1-\qi}(x_i)\leq r\leq r_i.\eqno\duediciotto$$
Suppose that the blow-up is not simple; then there exist some sequences of positive numbers $\tilde r_i\to 0$, $\tilde c_i\to\infty$, satisfying $\tilde c_i\vi^{1-\qi}(x_i)\leq \tilde r_i$ such that after passing to a subsequence
$$\bar w_i'(\tilde r_i)\geq 0.\eqno\duediciannove$$
It follows from \duediciotto\ and \duediciannove\ that $\tilde r_i\geq r_i$ and $\bar w_i$ has at least one critical point in the interval $[r_i,\tilde r_i]$. Let $\mu_i$ be the smallest critical point of $\bar w_i$ in this interval. It is clear that
$$\tilde r_i\geq\mu_i\geq r_i\quad\hbox{and}\quad\lim_{i\to\infty}\mu_i=0.$$
Consider now
$$\xi_i(x)=\mu_i^{1\over\qi-1}\vi(\mu_i x+x_i).$$
Set $T_i=x_i^n/\mu_i$ and $T=\lim_iT_i$. Then we have that $\xi_i$ satisfies the following
$$\left\{\eqalign{&-\Delta\xi_i=0,\quad\xi_i>0,\qquad\quad\hbox{in}\ B^{-T_i}_{1/\mu_i},\cr
                  &-{\partial \xi_i\over\partial x^n}=(n-2)f_i^{\ti}\xi_i^{\qi},\quad\ \ \hbox{on}\ \partial B^{-T_i}_{1/\mu_i}\cap \{x^n=-T_i\},\cr
                  &|x|^{1\over \qi-1}\xi_i(x)\leq C_10,\ \qquad\qquad\hbox{in}\ B^{-T_i}_{1/\mu_i},\cr
&\lim_{i\to\infty}\xi_i(0)=\infty\quad\hbox{and}\ 0\ \hbox{is a local maximum point of}\ \xi_i,\cr
     &r^{1\over\qi-1}\bar\xi_i(r)\quad\ \ \ \ \ \ \ \hbox{has negative derivative in}\quad C_10\xi_i(0)^{1-\qi}<r<1,\cr
   &\left.{d\over dr}\left(r^{1\over \qi-1}\bar\xi_i(r)\right)\right|_{r=1}=0.\cr}\right.\eqno\lunga$$
It is easy, arguing as we did before (e.g. see the proof of Lemma 2.9), to see that $\{\xi_i\}_i$ is locally bounded and then converges to some function $\xi$ satisfying
$$\left\{\eqalign{&-\Delta\xi=0,\quad\xi>0,\quad\quad\ \ \hbox{in}\ \erre^n_{-T},\cr
&-{\partial \xi\over\partial x^n}=(n-2)\xi^{n\over n-2},\quad \hbox{on}\ \partial\erre^n_{-T}.\cr}\right.$$
By the Liouville Theorem and the uniqueness result by Li and Zhu [\lizhuuno] of the appendix we deduce that $T=0$. Since $0$ is an isolated simple blow-up point, by Lemma 2.9 we have that
$$\xi_i(0)\xi_i(x)\ \mathop{\longrightarrow}\limits_{i\to\infty}\ a|x|^{2-n}+b=h(x)\quad\hbox{in}\quad C^2_{\rm loc}(B_1^+\privato\{0\})\eqno\dueventi$$
where $a>0$ and $b$ is some harmonic function satisfying
$$\left\{\eqalign{&-\Delta b=0,\quad\hbox{in}\ \erre^n_+,\cr
&{\partial b\over\partial x^n}=0,\quad\ \ \hbox{on}\ \partial\erre^n_+\privato\{0\}.\cr}\right.$$
By the Maximum Principle we see that $b\geq0$. Now, reflecting $b$ to be defined on all $\erre^n$ and denoting the resulting function by $\tilde b$, we deduce from the Liouville Theorem that $\tilde b$ is a constant and so $b$ is a constant. Using the last equality in \lunga\ and \dueventi, we deduce easily that $a=b$. Hence $h(x)=a(|x|^{2-n}+1)$. Therefore by Corollary 6.4 in the appendix we have that
$$\lim_{r\to0}\int_{\Gamma_2(B_r^+)}B(r, x, h,\nabla h)<0\eqno\dueventuno$$
where $B$ is given by
$$B(x,r,h,\nabla h)={n-2\over2}{\partial h\over\partial\nu}h+\mez r\left({\partial h\over\partial\nu}\right)^2-\mez r|\nabla_{\rm tan}h|^2\eqno\defdibi$$
where $\nabla_{\rm tan}h$ is the tangent component of $\nabla h$. From another part, using Lemma 6.3 in the appendix, Lemma 2.6, and Lemma 2.9, we deduce
$$\int_{\Gamma_2(B_r^+)}B(r, x, \xi_i,\nabla \xi_i)\geq O\left(\vi^{-2}(x_i)\right)\ti+O\left(\vi^{-(\qi+1)}(x_i)\right).$$
Multiplying by $\xi_i^2(0)$ we derive that
$$\lim_{r\to0}\int_{\Gamma_2(B_r^+)}B(r, x, h,\nabla h)\geq0,$$
which is in contradiction with \dueventuno. Therefore our Proposition is proved.\hfill\fine
\goodbreak

\vskip2truecm
\nuovoparagrafo 3

\centerline{\bf 3.\quad Compactness results for manifolds of positive type}
\bigskip
We point out that if $q$ stays strictly below the critical exponent ${n\over n-2}$ and strictly above $1$, the compactness of solutions of $({\cal P}_q)$ is much easier matter since it follows directly from the nonexistence of positive solutions to the global equation which one arrives at after a rather standard blow-up argument. Namely we prove
\medskip\noindent
{\bf Theorem 3.1.}\quad \it Let $(M,g)$ be a smooth compact $n\!$-dimensional Riemannian manifold with boundary. Then for any $\delta_1>0$ there exists a constant $C=C(M,g,\delta_1)>0$ such that for all $u\in\bigcup_{1+\delta_1\leq q\leq{n\over n-2}-\delta_1}{\cal M}_q^+$ we have
$${1\over C}\leq u(x)\leq C,\quad\forall\, x\in M;\quad\|u\|_{C^2(M)}\leq C.$$
\medskip\noindent
{\bf Proof.}\quad\rm Suppose that the Theorem were false. Then, in view of the Harnack inequality (see Lemma 6.1 in the appendix) and standard elliptic estimates, we would find sequences $\{q_i\}_i$ and $\{u_i\}_i\subset{\cal M}_{\qi}$ satisfying
$$\lim_{i\to\infty}\qi=q\in\bigg]1,{n\over n-2}\bigg[\quad\hbox{and}\quad \lim_{i\to\infty}\max_M\ui=\infty.$$
Let $p_i$ be the maximum point of $\ui$; it follows from the Maximum Principle that $p_i\in\partial M$. Let $x$ be a geodesic normal coordinate system in a neighbourhood of $p_i$ given by $\exp_{p_i}^{-1}$. We write $\ui(x)$ for $\ui(\exp_{p_i}(x))$. We rescale $x$ by $y=\lambda_i x$ with $\lambda_i=\ui^{\qi-1}(p_i)\to\infty$ and define 
$$\hat v_i(y)=\lambda_i^{-{1\over \qi-1}}\ui(\lambda_i^{-1}y).$$
Clearly $\hat v_i(0)=1$ and $0\leq\hat v_i\leq 1$. Let $\delta>0$ be some small positive number independent of $i$. We write $g(x)=g_{ab}(x)\,dx^adx^b$ for $x\in\exp_{p_i}^{-1}(B_\delta(p_i)\cap M)$. Define
 $$g^{(i)}(y)=g_{ab}(\lambda_i^{-1}y)\, d y^ady^b.$$
 Then $\hat \vi$ satisfies
$$\left\{\eqalign{&-L_{ g^{(i)}}\hat\vi=0,\quad\hat\vi>0,\quad\hbox{in}\ \lambda_i\exp^{-1}(B_\delta(p_i)\cap M),\cr
&B_{g^{(i)}}\hat\vi=(n-2)\hat\vi^{\qi},\ \quad\quad\hbox{on}\  \lambda_i\exp^{-1}(B_\delta(p_i)\cap\partial M).\cr}\right.$$
Applying $L^p\!$-estimates and Schauder estimates, we know that, after passing to a subsequence and a possible rotation of coordinates , $\hat\vi$ converges to a limit $\hat v$ in $C^2\!$-norm on any compact subset of $\{y\in\erre^n:\ y^n\geq0\}$, where
$$\left\{\eqalign{&-\Delta\hat v=0,\quad\ \ \  \quad\qquad\hbox{in}\ \erre^n_+,\cr
&-{\partial\hat v\over\partial y^n}=(n-2)\hat v^q,\quad\hbox{on}\ \partial \erre^n_+,\cr
&\hat v(0)=1,\quad \hat v\geq 0.\cr}\right.\eqno\treuno$$
It follows from the Liouville Theorem by Hu [\hu] that \treuno\ has no solution. This is a contradiction, thus we have established Theorem 3.1.\hfill\fine\medskip
The compactness of solutions of $({\cal P}_q)$ is much more difficult to establish when allowing $q$ to be close to ${n\over n-2}$, since the corresponding global equation does have solutions. On the other hand, due to the Liouville Theorem and Liouville-type Theorem by Li-Zhu [\lizhuuno] on the half-space $\erre^n_+$, we have the following Proposition similar to Lemma 3.1 of [\schoenzhang] and Proposition 1.1 of~[\hanliuno]. 
\medskip\noindent
{\bf Proposition 3.2.}\quad \it Let $(M,g)$ be a smooth compact $n\!$-dimensional Riemannian ma\-ni\-fold with boundary. For any $R\geq1$, $0<\eps<1$, there exist positive constants $\delta_0=\delta_0(M,g,R,\eps)$, $c_0=c_0(M,g,R,\eps)$, and $c_1=c_1(M,g,R,\eps)$ such that for all $u$ in
$$\bigcup_{{n\over n-2}-\delta_0\leq q\leq{n\over n-2}}{\cal M}_q^+$$
with $\max_Mu\geq c_0$, there exists ${\cal S}=\{p_1,\dots,p_N\}\subset \partial M$ with $N\geq 1$ such that 
\smallskip\noindent
\item{\rm(i)} each $p_i$ is a local maximum point of $u$ in $M$ and
$$\overline{B_{\bar r_i}(p_i)}\cap\overline{B_{\bar r_j}(p_j)}=\emptyset,\quad\hbox{for}\ i\not=j,$$
where $\bar r_i=Ru^{1-q}(p_i)$ and $\overline{B_{\bar r_i}(p_i)}$ denotes the geodesic ball in $(M,g)$ of radius $\bar r_i$ and centered at $p_i$;
\smallskip\noindent
\item{\rm(ii)} $$\left\|u^{-1}(p_i)u(\exp_{p_i}(yu^{1-q}(p_i))-\left({1\over (1+x^n)^2+|x'|^2}\right)^{n-2\over 2}\right\|_{C^2(B_{2R}^M(0))}<\eps$$
where 
$$B_{2R}^M(0)=\{y\in T_{p_i}M:\ |y|\leq 2R,\ u^{1-q}(p_i)y\in\exp^{-1}_{p_i}(B_\delta(p_i))\},$$
 $y=(y',y^n)\in\erre^n$;
\smallskip\noindent
\item{\rm(iii)} $d^{1\over q-1}(p_j,p_i)u(p_j)\geq c_0$, for $j>i$, while $d(p,{\cal S})^{{1\over q-1}}u(p)\leq c_1$, $\forall \, p\in M$, where $d(\cdot,\cdot)$ denotes the distance function in metric $g$.
\medskip\noindent\rm
The proof of Proposition 3.2 will follow from the following Lemma.
\medskip\noindent
{\bf Lemma 3.3.}\quad\it Let $(M,g)$ be a smooth compact $n\!$-dimensional Riemannian manifold. Given $R\leq1$ and $\eps<1$, there exist positive constants $\delta_0=\delta_0(M,g,R,\eps)$ and $C_0=C_0 (M,g,R,\eps)$ such that, for any compact $K\subset M$ and any $u\in\bigcup_{{n\over n-2}-\delta_0\leq q\leq{n\over n-2}}{\cal M}_q$ with $\max_{p\in\overline{M\privato K}}d^{1\over q-1}(p, K)u(p)\leq C_0$, we have that there exists $p_0\in M\privato K$ which is a local maximum point of $u$ in $M$ such that $p_0\in\partial M$ and
$$\left\|u^{-1}(p_0)u(\exp_{p_0}(yu^{1-q}(p_0))-\left({1\over (1+x^n)^2+|x'|^2}\right)^{n-2\over 2}\right\|_{C^2(B_{2R}^M(0))}<\eps$$
where $B_{2R}^M(0)$ is as in Proposition 3.2, $d(p,K)$ denotes the distance of $p$ to $K$, with $d(p,K)=1$ if $K=\emptyset$.
\medskip\noindent\rm
{\bf Proof.}\quad Suppose the contrary, then there exist compacta $K_i\subset M$, ${n\over n-2}-{1\over i}\leq\qi\leq{n\over n-2}$, and solutions $\ui$ of $({\cal P}_{\qi})$ such that 
$$\max_{p\in\overline{M\privato K_i}}d^{1\over \qi-1}(p,K_i)\ui(p)\geq i.$$
It is easy to deduce from the Hopf Lemma that $\ui>0$ in $M$. Let $\hat p_i\in\overline{M\privato K_i}$ be such that 
$$d^{1\over \qi-1}(\hat p_i,K_i)\ui(\hat p_i)=\max_{p\in\overline{M\privato K_i}}d^{1\over \qi-1}(p,K_i)\ui(p).$$
Let $x$ be a geodesic normal coordinate system in a neighbourhood of $\hat p_i$ given by $\exp^{-1}_{\hat p_i}$. We write $\ui(x)$ for $\ui(\exp_{\hat p_i}(x))$ and denote $\lambda_i=\ui^{\qi-1}(p_i)$. We rescale $x$ by $y=\lambda_i x$ and define $\hat\vi(y)=\lambda_i^{-{1\over\qi-1}}\ui(\lambda_iy)$. By standard blow-up arguments and the Liouville Theorem, one can prove that $d(\hat p_i,\partial M)\to0$. Fix some small positive constant $\delta>0$ independent of $i$ such that $\partial M\cap B_\delta(\hat p_i)\not=\emptyset$. We may assume without loss of generality, by taking $\delta$ smaller, that $\exp_{\hat p_i}^{-1}(\partial M)\cap B_{\delta}(0)$ has only one connected component, and may arrange to let the closest point on $\exp_{\hat p_i}^{-1}(\partial M)\cap B_{\delta}(0)$ to $0$ to be at $(0,\dots,0,-t_i)$ and $$\exp_{\hat p_i}^{-1}(\partial M)\cap B_{\delta}(0)=\partial \erre^n_+\cap B_\delta^M(0)$$ is a graph over $(x^1,\dots,x^{n-1})$ with horizontal tangent plane at $(0,\dots,-t_i)$ and uniformly bounded second derivatives. In $\exp^{-1}_{\hat p_i}(B_\delta(\hat p_i))$ we write $g(x)=g_{ab}(x)\,dx^adx^b$. Define 
$$g^{(i)}(y)=g_{ab}(\lambda_i^{-1}y)\, d y^ady^b.$$
 Then $\hat \vi$ satisfies
$$\left\{\eqalign{&-L_{ g^{(i)}}\hat\vi=0,\quad\hat\vi>0,\cr
&B_{g^{(i)}}\hat\vi=(n-2)\hat\vi^{\qi}.\cr}\right.$$
Note that $\lambda_id(\hat p_i, K_i)\to\infty$ and, for $|y|\leq\quarto\lambda_id(\hat p_i, K_i)$ with $x=\lambda_i^{-1}y\in\exp_{\hat p_i}^{-1}(B_\delta(\hat p_i))$, we have
$$d(x,K_i)\geq\mez d(\hat p_i, K_i),$$
and therefore
$$\left(\mez d(\hat p_i, K_i)\right)^{1\over \qi-1}\ui(x)\leq d(x,K_i)^{1\over \qi-1}\ui(x)\leq d(\hat p_i, K_i)^{1\over \qi-1}\ui(\hat p_i)$$
which implies , for all $|y|\leq \quarto\lambda_id(\hat p_i, K_i)$ with $\lambda_i^{-1}y\in\exp^{-1}_{\hat p_i}(B_\delta(\hat p_i))$, that 
$$\hat v_i(y)\leq 2^{1\over\qi-1}.$$
Standard elliptic theories imply that there exists a subsequence, still denoted by $\hat\vi$, such that, for $T=\lim_i\lambda_id(\hat p_i,\partial M)\in[0,+\infty]$, $\hat \vi$ converges to a limit $\hat v$ in $C^2\!$-norm on any compact set of $\{y=(y^1,\dots, y^n)\in\erre^n:\ y^n\geq-T\}$, where $\hat v>0$ satisfies
$$\left\{\eqalign{&-\Delta \hat v=0,\ \quad\qquad\qquad\quad\hbox{in}\ \{y^n>-T\},\cr
&-{\partial \hat v\over\partial y^n}=(n-2)\hat v^{n\over n-2},\quad\hbox{on}\ \{y^n=-T\},\quad\hbox{if}\ T<+\infty.\cr}\right.$$
It follows from the Liouville Theorem that $T<+\infty$, and, from the Liouville-type Theorem of Li-Zhu [\lizhuuno], that
$$\hat v(x',x^n)=\left({1\over (1+(x^n-T))^2+|x'-x_0'|^2}\right)^{n-2\over2}.$$
Set $\hat y=(\hat y',-T)$. It follows from the explicit form of $\hat\vi$ that there exist $\yi\to\hat y$ which are local maximum points of $\hat\vi$ such that $\hat\vi(\yi)\to \lambda^{n-2\over 2}=\max\hat v$.\par
Define $p_i=\exp_{\hat p_i}(\lambda_i^{-1}\yi)$, then $p_i\in M\privato K_i$ is a local maximum point of $\ui$, and if we repeat the scaling with $p_i$ replacing $\hat p_i$, we still obtain a new limit $v$. Due to our choice, $v(0)=1$ is a local maximum, so $T=0$ and
$$\left\|\ui^{-1}(p_i)\ui(\exp_{p_i}(y\ui^{1-\qi}(p_i))-\left({1\over (1+x^n)^2+|x'|^2}\right)^{n-2\over 2}\right\|_{C^2(B_{2R}^M(0))}<\eps$$
which leads to a contradiction.\hfill\fine
\medskip\noindent
{\bf Proof of Proposition 3.2.}\quad First we apply Lemma 3.3 by taking $K=\emptyset$ and $d(p,K)\equiv1$ to obtain $p_1\in\partial M$ which is a maximum point of $u$ and (i) of Lemma 3.3 holds. If
$$\max_{p\in M\privato K_1}d^{1\over q-1}(p,K_1)u(p)\leq C_0,$$
where $K_1=\overline{B_{\bar r_1}(p_1)}$, we stop. Otherwise we apply again Lemma 3.3 to obtain $p_2\in\partial M$. It is clear that we have $\overline{B_{\bar r_1}(p_1)}\cap\overline{B_{\bar r_2}(p_2)}=\emptyset$ by taking $\eps$ small from the beginning. We continue the process. Since there exists $a(n)>0$ such that $\int_{B_{\bar r_i}(p_i)}\ui^{\qi+1}\geq a(n)$, our process will stop after a finite number of steps. Thus we obtain ${\cal S}=\{p_1,\dots,p_N\}\subset\partial M$ as in (ii) and 
$$d^{1\over q-1}(p,{\cal S})u(p)\leq C_0,$$
for any $p\in M\privato{\cal S}$. Clearly, we have that item (iii) holds.\hfill\fine
\medskip\noindent
Though Proposition 3.2 states that $u$ is very well approximated in strong norms by standard bubbles in disjoint balls $B_{\bar r_1}(p_1),\dots,B_{\bar r_N}(p_N)$, it is far from the compactness result we wish to prove. Interactions between all these bubbles have to be analyzed to rule out the possibility of blowing-ups.\par
The next Proposition rules out possible accumulations of these bubbles, and this implies that only isolated blow-up points may occur to a blowing-up sequence of solutions.
\medskip\noindent
{\bf Proposition 3.4.}\quad \it Let $(M,g)$ be a smooth compact $n\!$-dimensional locally conformally flat Riemannian manifold with umbilic boundary. For suitably large $R$ and small $\eps>0$, there exist $\delta_1=\delta_1(M,g,R,\eps)$ and $d=d(M,g,R,\eps)$ such that for all $u$ in 
$$\bigcup_{{n\over n-2}-\delta_1\leq q\leq{n\over n-2}}{\cal M}_q^+$$
with $\max_Mu\geq C_0$, we have
$$\min\{d(p_i,p_j):\ i\not=j,\ 1\leq i,j\leq N\}\geq d$$
where $C_0,p_1,\dots,p_N$ are given by Proposition 3.2.
\medskip\noindent
{\bf Proof.}\quad \rm By contradiction, suppose that the conclusion does not hold, then there exist sequences ${n\over n-2}-{1\over i}\leq \qi\leq{n\over n-2}$, $\ui\in{\cal M}_{\qi}$ such that $\min\{d(p_{i,a},p_{i,b}),\ 1\leq a,b\leq N\}\to0$ as $i\to+\infty$ where $p_{i,1},\dots,p_{i,N}$ are the points given by Proposition 3.2. Notice that when we apply Proposition 3.2 to determine these points, we fix some large constant $R$, and then some small constant $\eps>0$ (which may depend on $R$), and in all the arguments $i$ will be large (which may depend on $R$ and $\eps$). Let
$$d_i=d(p_{i,1},p_{i,2})=\min_{a\not=b}d(p_{i,a},p_{i,b})$$
and 
$$p_0=\lim_{i\to+\infty}p_{i,1}=\lim_{i\to+\infty}p_{i,2}\in\partial M.$$
Since $M$ is locally conformally flat with umbilic boundary, one can find a diffeomorphism 
$$\Phi:\quad B_2^+\,\longrightarrow\, B_\delta(p_0),\quad \Phi(0)=p_0\eqno\diffeo$$
with $\Phi^{\star}g=f^{4\over n-2}g_0$ where $g_0$ is the flat metric in $B_2^+$ and $f\in C^2(\overline{B_2^+})$ is some positive function. It follows from the conformal invariance of $L_g$ and $B_g$ that, for $v_i=f\ui\circ\Phi$,
$$\left \{\eqalign {&-\Delta \vi=0,\quad \vi>0,\quad\quad \hbox{in}\ B_2^+,\cr
&{\partial \vi\over\partial x^n}=-(n-2)f^{\ti}\vi^{\qi},\quad\,\hbox{on}\ \Gamma_1(B_2^+).\cr}\right.\eqno\diffeoprimo$$
We can assume without loss of generality that $x_{i,a}=\Phi^{-1}(p _{i,a})$ are local maxima of $\vi$, so it is easy to see that
$$\eqalignno{&\vi(x_{i,a})\,\longrightarrow\,+\infty,&\auno\cr
&d\bigg(x,\bigcup_a\{x_{i,a}\}\bigg)^{1\over\qi-1}\vi(x)\leq c_1,\quad\forall\, x\in B_1^+,&\adue\cr
&0<\sigma_i:=|x_{i,1}-x_{i,2}|\,\longrightarrow\,0,&\cr
&\sigma_i^{1\over \qi-1}\vi(x_i,x)\geq {R^{n-2\over 2}\over c_2}\quad\hbox{for}\ a=1,2,&\atre\cr}$$
where $c_1,c_2>0$ are some constants independent of $i,\eps,R$. Without loss of generality, we assume that $x_{i,1}=(0,\dots,x_{i,1}^n)$. Consider 
$$w_i(y)=\sigma_i^{1\over \qi-1}\vi(x_{i,1}+\sigma_iy)$$
and set, for $x_{i,a}\in\overline{B_1^+}$, $y_{i,a}={x_{i,a}-x_{i,1}\over\sigma_i}$ and $T_i={1\over \sigma_i}x^n_{i,a}$. Clearly, $w_i$ satisfies
$$\left\{\eqalign{&-\Delta w_i(y)=0,\quad w_i>0,\qquad\ \qquad\quad\hbox{in}\ \left\{|y|<{1\over\sigma_i},\ y^n>-T_i\right\},\cr
&{\partial w_i\over\partial y^n}=-(n-2)f^{\ti}(x_{i,1}+\sigma_iy)w_i^{\qi},\quad\!\hbox{on}\ \left\{|y|<{1\over\sigma_i},\ y^n=-T_i\right\}.\cr}\right.\eqno\tretre$$
It follows that
$$|y_{i,a}-y_{i,b}|\geq 1,\quad\forall\, a\not=b,\quad y_{i,1}=0,\quad|y_{i,2}|=1.\eqno\trequattro$$
After passing to a subsequence, we have
$$\bar y=\lim_{i\to+\infty}y_{i,2},\quad |\bar y|=1.$$
It follows easily from \auno, \adue, and \atre\ that
$$\left\{\eqalign{&w_i(0)\geq c_0'\quad w_i(y_{i,2})\geq c_0',\cr
&\hbox{each}\ y_{i,a}\ \hbox{is a local maximum point of}\ w_i,\cr
&\min_a|y-y_{i,a}|^{1\over\qi-1}w_i(y)\leq c_1,\cr
&|y|\leq{1\over 2\sigma_i},\quad y^n\geq -T_i\cr}\right.$$
where $c_0'>0$ is independent of $i$. At this point we need the following Lemma which is a direct consequence of Lemma 2.4.
\medskip\noindent
{\bf Lemma 3.5.}\quad\it If along some subsequence both $\{y_{i,a_i}\}$ and $w_i(y_{i,a_i})$ remain bounded, then along the same subsequence 
$$\limsup_{i\to+\infty}\max_{B^{-T_i}_{1/4}(y_{i,a_i})}w_i<\infty,$$
where $B^{-T_i}_{1/4}(y_{i,a_i})=\{y:\ |y-y_{i,a_i}|<1/4,\ y^n>-T_i\}$.
\medskip\noindent\rm
Due to~Proposition 2.11 and Lemma~3.5, all the points $y_{i,a_i}$ are either regular points of $w_i$ or isolated simple blow-up points. We deduce, using Lemma 2.9, Lemma 3.5, \tretre, and \trequattro \ that
$$w_i(0)\,\longrightarrow\,+\infty,\quad w_i(y_{i,2})\,\longrightarrow\,+\infty.$$
It follows that $\{0\},\ \{y_{i,2}\to\bar y\}$ are both isolated simple blow-up points. Let $\tilde w_i=w_i(0)w_i$. It follows from Lemma 2.9 that there exists $\widetilde{\cal S}_1$ such that $\{0,\bar y\}\subset \widetilde{\cal S}_1\subset{\cal S}$,
$$\min\{|x-y|:\ x,y\in\widetilde{\cal S}_1,\ x\not= y\}\geq1,$$
and
$$w_i(0)w_i\ \mathop{\longrightarrow}\limits_{i\to\infty}\ h\quad\hbox{in}\ C^2_{\rm loc}(\erre^n_{-T}\privato\widetilde{\cal S}_1)$$
where $h$ satisfies
$$\left\{\eqalign{&\Delta h=0,\quad \ \,\hbox{in}\ \erre^n_{-T}\privato\widetilde{\cal S}_1,\cr
                  &{\partial h\over\partial y^n}=0,\quad\hbox{on}\ \partial\erre^n_{-T}\privato\widetilde{\cal S}_1.\cr}\right.$$
Making an even extension of $h$ across the hyperplane $\{y^n=-T\}$, we obtain $\tilde h$ satisfying $\Delta\tilde h=0$ on $\erre^n\privato\widetilde{\cal S}_1$. Using B\"ocher's Theorem, the fact that $\{0,\bar y\}\subset\widetilde{\cal S}_1$, and the Ma\-xi\-mum Principle, we obtain some nonnegative function $b(y)$ and some positive constants $a_1,a_2>0$ such that
$$\left\{\eqalign{&b(y)\geq 0,\ \ \,\quad y\in \erre^n\privato\{\widetilde{\cal S}_1\privato\{0,\bar y\}\},\cr
&\Delta b(y)=0,\quad  y\in \erre^n\privato\{\widetilde{\cal S}_1\privato\{0,\bar y\}\},\cr
&{\partial b\over\partial\nu}=0,\quad\ \ \ \,  \hbox{on}\ \partial \erre^n_+\privato\{\widetilde{\cal S}_1\privato\{0,\bar y\}\},\cr}\right.$$
and $h(y)=a_1|x|^{2-n}+a_2|x-\bar y|^{2-n}+b$, $y\in\erre^n\privato\widetilde{\cal S}_1$. Therefore there exists $A>0$ such that
$$h(y)=a_1|y|^{2-n}+A+O(|y|)$$
for $y$ close to zero. Using Lemma 6.3 and Corollary 6.4 in the appendix, we obtain a contradiction as in Proposition 2.11. The proof of our Proposition is thereby complete.\hfill\fine
\medskip\noindent
{\bf Proof of Theorem 1.5.}\quad Let $f_1$ be an eigenfunction of problem $(E_1)$ associated to $\lambda_1(L)$. Taking if necessary $|f_1|$, we can assume $f_1\geq0$. By the Maximum Principle $f_1>0$ in $\mathop{M}\limits^{\circ}$ and by the Hopf Maximum Principle $f_1>0$ on $\partial M$. Thus $f_1>0$ in~$M$. Consider the metric $g_1=f_1^{4\over n-2}g$. Then $R_{g_1}>0$ and $h_{g_1}\equiv0$. We will work with $g_1$ instead of $g$. For simplicity of notation, we still denote it as $g$. Then we can assume $R_g>0$ and $h_g\equiv0$ without loss of generality, so that $B_g=\partial/\partial\nu$.\par
In view of $L^p\!$-estimates, Schauder estimates, and Lemma 6.1, we only need to esta\-blish the $L^\infty\!$-bound of $u$. Arguing by contradiction, suppose there exist sequences $\qi={n\over n-2}-\ti$, $\ti\geq0$, $\ti\to0$, and $\ui\in{\cal M}_{\qi}$ such that
$$\max_M\ui\ \mathop{\longrightarrow}\limits_{i\to\infty}\ \infty.$$
It follows from Proposition 2.10, Theorem 3.1, and Proposition 3.4 that, after passing to a subsequence, $\{\ui\}_i$ has $N$ $(1\leq N<\infty)$ isolated simple blow-up points denoted by $\{p^1,\dots,p^N\}$. Let  $\{p_i^1,\dots,p_i^N\}$ denote the local maximum points as in Definition 2.1. It follows from Proposition 2.5 that
$$\ui(p_i^1)\ui\ \mathop{\longrightarrow}\limits_{i\to\infty}\ h\quad\hbox{in}\ C^2_{\rm loc}(M\privato\{p^1,\dots,p^N\}).$$
Using Proposition 2.5 and subtracting to the function $h$ the contribution of all the poles $\{p^1,\dots,p^N\}\subset\partial M$, we obtain
$$\ui(p_i^1)\ui\ \mathop{\longrightarrow}\limits_{i\to\infty}\ \sum_{\ell=1}^Na_\ell G(\cdot, p^\ell)+\tilde b\quad\hbox{in}\ C^2_{\rm loc}(M\privato\{p^1,\dots,p^N\})$$
where $a_\ell>0$, $G(\cdot, p^\ell)$ is as in \green, and $\tilde b$ satisfies
$$\left\{\eqalign{&L_g\tilde b=0,\quad\hbox{in}\ M,\cr
&B_g\tilde b=0,\quad\hbox{on}\ \partial M.\cr}\right.$$
Since $\lambda_1(L)>0$ we deduce that $\tilde b=0$ and $G(\cdot,p^\ell)>0$ (recall that we have chosen $g$ such that $R_g>0$ and $h_g\equiv0$). Since $M$ is compact and locally conformally flat with umbilic boundary, for every $p^\ell$ there exist $\rho>0$ uniform and $g_2=f_2^{4\over n-2}g$, for $f_2\in C^2(\overline{B_{2\rho}(p^\ell)})$, such that $g_2$ is Euclidean in a neighbourhood of $p^\ell$ and $h_g=0$ on $\partial M\cap B_\rho(p^\ell)$. It is standard to see that the Green's function $\widehat G(x, p^\ell)$ of $g_2$ has the following expansion near $p^\ell$ in geodesic normal coordinates
$$ \widehat G(x, p^\ell)=|x|^{2-n}+A+O(|x|).$$
It follows then from the Positive Mass Theorem by Schoen and Yau [\schoenyau] as it was extended to locally conformally flat manifolds with umbilic boundary by Escobar [\escobartre] that $A\geq 0$ with equality if and only if $(M,g)$ is conformally equivalent to the standard ball. Let $\vi$ be as in \diffeoprimo. Recall that $\Phi(p^1)=0$, so we can deduce that $x_i\to0$ is an isolated simple blow-up point of $\{\vi\}_i$ and
$$\vi(x_i)\vi\ \mathop{\longrightarrow}\limits_{i\to\infty}\ \tilde h\quad\hbox {in}\ C^2_{\rm loc}(\overline{B_1^+}\privato\{0\})$$
where $\tilde h(x)=|x|^{2-n}+\widetilde A+O(|x|)$ for some $\widetilde A>0$. Applying Lemma 6.3 and Corollary 6.4 of the appendix, we reach as usual a contradiction. The Theorem is then proved.\hfill\fine

  \vskip2truecm 
\goodbreak

\centerline{\bf 4.\quad Existence results for manifolds of positive type}
\bigskip
\rm In this section we prove the existence part of Theorem 1.4, using the compactness results of the previous section and the Leray-Schauder degree theory.\par
We assume $R_g>0$ and $h_g\equiv0$ without loss of generality (see the beginning of the proof of Theorem 1.5) so that $B_g=\partial/\partial\nu$. For $1\leq q\leq{n\over n-2}$, consider the problem
$$\left\{\eqalign{&L_gu=0,\quad\, \quad\qquad\hbox{in}\ \mathop{M}\limits^\circ,\cr
&{\partial u\over\partial\nu}=v,\, \ \quad\quad\qquad\hbox{on}\ \partial M,\cr}\right.\eqno{({\cal P}_v)}$$
which defines an operator\nuovoparagrafo 4
$$\eqalign{T:\quad C^{2,\alpha}(M)^+&\ \longrightarrow\ C^{2,\alpha}(M)\cr
                               v      &\ \longmapsto\ Tv=u\cr}$$
where $C^{2,\alpha}(M)^+:=\{u\in C^{2,\alpha}(M):\ u>0$ in $M\}$, $0<\alpha<1$ and $Tv$ is the unique solution of problem $({\cal P}_v)$. Set
$$E(v):=\int_M(-L_gv)v+\int_{\partial M}(B_gv)v=\int_M|\nabla_gv|^2+{n-2\over 4(n-1)}\int_MR_g v^2$$
and consider the problem
$$\left\{\eqalign{&-L_gv=0,\quad v>0,\ \ \quad\hbox{in}\ \mathop{M}\limits^{\circ},\cr
&B_gv=(n-2)E(v)v^q,\quad\hbox{on}\ \partial M.\cr}\right.\eqno\quattrouno$$
\medskip\noindent
We have the following Lemma
\medskip\noindent
{\bf Lemma 4.1.}\quad \it There exists some positive constant $C=C(M,g)$ such that, for all $1\leq q\leq {n\over n-2}$ and $v$ satisfying \quattrouno, we have
$${1\over C}<v<C,\quad\hbox{in}\ M.\eqno\quattrodue$$
\medskip\noindent
{\bf Proof.}\quad\rm First of all, notice that, in view of the Harnack inequality and Lemma 6.1, it is enough to prove the upper bound. Multiplying \quattrouno\ by $v$ and integrating by parts, we obtain
$$(n-2)E(v)\int_{\partial M}v^{q+1}=\int_M|\nabla_gv|^2+{n-2\over 4(n-1)}\int_M R_gv^2\eqno\quattrotre$$
which yields $E(v)>0$. It is easy to check that $u=E(v)^{1\over q-1}v>0$ satisfies
$$\left\{\eqalign{&-L_gu=0,\quad u>0,\quad\hbox{in}\ \interno{M},\cr
&B_gu=(n-2)u^q,\quad\ \quad\hbox{on}\ \partial M.\cr}\right.$$
It follows from Theorem 3.1 and Proposition 3.4 that there exists $\delta_0>0$ such that for $1+\delta_0\leq q\leq{n\over n-2}$
$${1\over c_1}\leq E(v)^{1\over q-1}v\leq c_1\eqno\quattroquattro$$
for some positive constant $c_1$. From \quattrotre\ we know that $(n-2)E(v)\int_{\partial M}v^{q+1}=E(v)$, so that 
$$\int_{\partial M}v^{q+1}={1\over n-2}.\eqno\quattrocinque$$
Next \quattroquattro\ and \quattrocinque\ yield 
$${1\over c_2}\leq E(v)\leq c_2\eqno\quattrosei$$
for some positive $c_2$. Then \quattroquattro\ and \quattrosei\ give \quattrodue\ for $1+\delta_0\leq q\leq {n\over n-2}$. For $1\leq q\leq 1+\delta_0$ we apply Lemma 6.5 to obtain $E(v)\leq c_3$ for a positive constant $c_3$ and then standard elliptic estimates to obtain the upper bound for $v$.\hfill\fine
\medskip\noindent
For $0<\alpha<1$, $1\leq q\leq{n\over n-2}$, we define a map 
$$\eqalign{F_q:\quad C^{2,\alpha}(M)^+&\ \longrightarrow\ C^{2,\alpha}(M)\cr
                               v      &\ \longmapsto\ F_qv=v-T(E(v)v^q).\cr}$$
For $\Lambda>1$, let
$$D_\Lambda=\left\{v\in C^{2,\alpha}(M),\ \|v\|_{C^{2,\alpha}(M)}<\Lambda,\ \min_M v>{1\over \Lambda}\right\}.\eqno\domain$$ 
Let us notice that $F_q$ is a Fredholm operator and $0\not\in F_q(\partial D_\Lambda)$ thanks to Lemma 4.1. Consequently, by the homotopy invariance of the Leray-Schauder degree (see [\nirenberg] for a comprehensive introduction to Leray-Schauder degree and its properties), we have
$$\deg (F_q,D_\Lambda,0)=\deg (F_1,D_\Lambda,0),\quad \forall\,1\leq q\leq{n\over n-2}.$$
It is easy to see that $F_1(v)=0$ if and only if $E(v)=\lambda_1(B)$ and $v=\root\of{\lambda_1(B)}f_2$, where $f_2$ is an eigenfunction of $(E_2)$ associated to $\lambda_1(B)$. Let $\bar v=\root\of{\lambda_1(B)}f_2$.
\medskip\noindent
{\bf Lemma 4.2.}\quad \it $F_1'(\bar v)$ is invertible with exactly one simple negative eigenvalue. Therefore $\deg(F_1,D_\Lambda,0)=-1$.
\medskip\noindent\rm
{\bf Proof.} This can be proved by quite standard arguments, one can follow, up to minor modifications, the derivation of similar results in [\hanliuno, pp. 528-529]. We omit the proof.~\hfill\fine
\medskip\noindent
For $s\in[0,1]$, let us consider the homotopy
$$\eqalign{G_s:\quad C^{2,\alpha}(M)^+&\ \longrightarrow\ C^{2,\alpha}(M)\cr
                               v      &\ \longmapsto\ G_s(v)=v-T_{n\over n-2}\left([(n-2)s+(1-s)E(v)]v^{n\over n-2}\right).\cr}$$  
Arguing as in Lemma 4.1, one easily deduces
\medskip\noindent
{\bf Lemma 4.3.}\quad\it There exists $\overline\Lambda>2$ depending only on $(M,g)$ such that
$$G_s(u)\not=0\quad\forall\, \Lambda\geq\overline\Lambda,\quad \forall\, 0\leq s\leq 1,\quad \forall\, u\in\partial D_\Lambda.$$
\medskip\noindent\rm
{\bf Proof of Theorem 1.4 completed.}\quad  Using Lemma 4.3 and the homotopy invariance of the Leray-Schauder degree, we have for all $\Lambda\geq\overline\Lambda$, 
$$\deg(G_1,D_\Lambda,0)=\deg(G_0,D_\Lambda,0).$$
Observing that
$$\eqalign{&G_1(u)=u-T_{n\over n-2}\left((n-2)u^{n\over n-2}\right),\cr
&G_0(u)=F_{n\over n-2}(u)\cr}$$
and using Lemma 4.2, we have that for $\Lambda$ sufficiently large
$$\deg(G_1,D_\Lambda,0)=-1,$$
which, in particular, implies that ${\cal M}\cap D_\Lambda\not=\emptyset$. We have thus completed the proof of the existence part of Theorem 1.4.\hfill\fine

\vskip2truecm

\centerline{\bf 5.\quad Compactness and existence results for manifolds of negative type}
\bigskip
In this section we establish Theorem 1.6. Let $f_2$ be a positive eigenfunction of $(E_2)$ corresponding to $\lambda_1(B)$ and set $g_2~=~f_2^{4\over n-2}g$. It follows that $R_{g_2}\equiv0$ and $h_{g_2}<0$. We will work throughout this section with $g_2$ instead of $g$ and we still denote it by $g$. \par
We first prove compactness part in Theorem 1.6. Due to the Harnack inequality, Lemma~6.1, elliptic estimates, and Schauder estimates, we need only to establish the $L^\infty\!$-bound. We use a contradiction argument. Suppose the contrary, that there exist sequences $\{\qi\}_i$, $\{\ui\}_i\in{\cal M}^-_{\qi}$ satisfying
$$\qi\ \mathop{\longrightarrow}\limits_{i\to\infty}\ q_0\in\left]1,{n\over n-2}\right]\quad\hbox{and}\quad  \lim_{i\to\infty}\max_M\ui=+\infty.$$
Let $x_i\in\partial M$ such that $\ui(x_i)=\max_M\ui\to+\infty$. Let $y^1,\dots,y^n$ be the geodesic normal coordinates given by some exponential map, with $\partial/\partial y^n=-\nu$ at $x_i$. Consider 
$$\tilde \ui(z)=\ui^{-1}(x_i)\ui\left(\exp_{x_i}(\ui^{1-\qi}(x_i)z)\right).$$
Reasoning as in Theorem 3.1, we obtain that $\tilde\ui$ converges in $C^2_{\rm loc}\!$-norm to some $\tilde u$ satisfying
$$\left\{\eqalign{&-\Delta\tilde u=0,\quad \tilde u>0,\quad\hbox{in}\ \erre^n_+,\cr
&{\partial \tilde u\over \partial z^n}=(n-2)u^{q_0},\ \;\quad\hbox{on}\ \partial\erre^n_+,\cr}\right.\eqno\cinqueuno$$
with $\tilde u(0)=1$, $0<\tilde u\leq 1$ on $\erre^n_+$. Using the Liouville-type Theorem of Lou-Zhu [\louzhu], we obtain that \cinqueuno\ has no solution satisfying $\tilde u(0)=1$ and $0<\tilde u\leq 1$. \par
We prove now the existence part of Theorem 1.6. Let 
$$E(u,v)=\int_M\nabla_g u\cdot\nabla_g v+{n-2\over 2}\int_{\partial M}h_guv$$
and $E(u)=E(u,u)$. Let us observe that one can choose $f_2$ such that $E(f_2)=-1$. Consider for, $1\leq q\leq{n\over n-2}$,
$$\left\{\eqalign{&\Delta_gv=0,\quad v>0,\quad\hbox{in}\ \interno M,\cr
&B_gv=E(v)v^q,\quad\;\quad\hbox{on}\ \partial M.\cr}\right.\eqno\cinquedue$$
Arguing as in Lemma 4.1 and using Lemma 6.6 one can prove
\medskip\noindent
{\bf Lemma 5.1.}\quad \it There exists some constant $C=C(M,g)>0$ such that for $1\leq q\leq{n\over n-2}$ and $v$ satisfying \cinquedue\ we have 
$${1\over C}<v<C.$$
\medskip\noindent
\rm
Let $\lambda_1(B)<\lambda_2(B)<\dots$ denote all the eigenvalues of $(E_2)$. Pick some constant $A~\in~(-\lambda_2(B),-\lambda_1(B))$. For $0<\alpha<1$ and $1\leq q\leq{n\over n-2}$, we define
$$\widetilde T_A:\quad C^{2,\alpha}(M)^+\ \longrightarrow\ C^{2,\alpha}(M),$$
which associates to $v\in C^{2,\alpha}(M)^+$ the unique solution of
$$\left\{\eqalign{&L_gu=0,\,\ \ \quad\qquad\hbox{in}\ \interno M,\cr
&(B_g+A)u=v,\quad\hbox{on}\ \partial M\cr}\right.$$
and $F_q(v)=v-\widetilde T_A(E(v)v^q+Av)$. For $\Lambda>1$, let $D_\Lambda\subset C^{2,\alpha}(M)^+$ be given as in \domain. It follows from Lemma 5.1 that $0\not\in F_q(\partial D_\Lambda)$, for all $1\leq q\leq{n\over n-2}$. Consequently,
$$\deg(F_q, D_\Lambda,0)=\deg(F_1,D_\Lambda,0),\quad\forall\, 1\leq q\leq{n\over n-2}.$$
Arguing as we did in Lemma 4.2, we obtain
\medskip\noindent
{\bf Lemma 5.2.}\quad\it Suppose $\lambda_1(B)<0$ and $R_g\equiv0$. Then
$$\deg(F_q,D_\Lambda,0)=-1,\quad\forall\,1\leq q\leq{n\over n-2}.$$
\medskip\noindent\rm
Now we define for $1\leq q\leq{n\over n-2}$, $\widetilde T_q$ as follows
$$\eqalign{\widetilde T_q:\quad C^{2,\alpha}(M)^+&\ \longrightarrow\ C^{2,\alpha}(M)\cr
                               v      &\ \longmapsto\ \widetilde T_qv=u\cr}$$
where $u$ is the unique solution of
$$\left\{\eqalign{&\Delta_gu=0,\ \quad\qquad\qquad\qquad\qquad\quad\hbox{in}\ \interno M,\cr
&(B_g+A)u=-(n-2)v^q+Av\quad\hbox{on}\ \partial M.\cr}\right.$$
Since $0$ is not an eigenvalue of $B_g+A$, $\widetilde T_q$ is well defined. It follows from Schauder theory, see e.g. [\gilbarg], that $\widetilde T_q$ is compact. It follows from the compactness part of Theorem~1.6 that there exists $\Lambda>\!>1$ depending only on $(M,g)$ such that
$$\left\{u\in C^{2,\alpha}(M):\ \Big({\rm Id}-\widetilde T_{n\over n-2}\Big)u=0\right\}\subset D_\Lambda\quad\hbox{for every}\ \Lambda>\bar\Lambda.$$
\medskip\noindent
{\bf Lemma 5.3.}\quad \it Suppose that $\lambda_1(B)<0$ and $R_g\equiv 0$. Then for $\Lambda$ large enough, we have
$$\deg\left({\rm Id}-T_{n\over n-2},D_\Lambda,0\right)=\deg(F_1,D_\Lambda,0)=-1.$$
\medskip\noindent
{\bf Proof.}\quad\rm It follows from the homotopy invariance of the Leray-Schauder degree from one part and Lemma 5.2 from another part. The proof being standard, we omit it.\hfill\fine
\medskip\noindent
{\bf Proof of Theorem 1.6 completed.}\quad The existence part follows from Lemma 5.3 and standard degree theory. Thereby the proof of Theorem 1.6 is established.\hfill\fine

\vskip2truecm
\nuovoparagrafo 6

\centerline{\bf Appendix}
\bigskip\noindent
In this appendix, we present some results used in our arguments. First of all we state a Harnack inequality for second-order elliptic equations with Neumann boundary condition. For the proof one can see [\hanliuno, Lemma A.1].
\medskip\noindent
{\bf Lemma 6.1.}\quad\it Let $L$ be the operator 
$$Lu=\partial_i(a_{ij}(x)\partial_j u+b_i(x)u)+c_i(x)\partial_iu+d(x)u$$
and assume that for some constant $\Lambda>1$ the coefficients satisfy
$$\eqalignno{&\Lambda^{-1}|\xi|^2\leq a_{ij}(x)\xi_i\xi_j\leq \Lambda|\xi|^2,\quad\quad\forall x\in B_3^+\subset\erre^n,\qquad\xi\in\erre^n,&\seiuno\cr
&|b_i(x)|+|c_i(x)|+|d(x)|\leq\Lambda,\qquad\,\forall \, x\in B_3^+.&\seidue\cr}$$
If $|h(x)|\leq\Lambda$ for any $x\in B_3^+$ and $u\in C^2(B_3^+)\cap C^1(\overline{ B_3^+})$ satisfies
$$\left\{\eqalign{&-Lu=0,\quad u>0,\quad \hbox{in}\ B_3^+,\cr
&a_{nj}(x)\partial_ju=h(x)u,\quad\hbox{on}\ \Gamma_1(B_3^+),\cr}\right.$$
then there exists $C=C(n,\Lambda)>1$ such that
$$\max_{\overline{B_1^+}}u\leq C\min_{\overline{B_1^+}}u.$$
\medskip\noindent\rm
In the proofs of our results, we also used the following Maximum Principle.
\medskip\noindent
{\bf Theorem 6.2.}\quad \it Let $\Omega$ be a bounded domain in $\erre^n$ and let $\partial \Omega=\Gamma\cup\Sigma$, $V\in L^{\infty}(\Omega)$, and $h\in L^\infty(\Sigma)$. Suppose $\psi\in C^2(\Omega)\cap C^1(\bar\Omega)$, $\psi>0$ in $\bar\Omega$ satisfies
$$\left\{\eqalign{&\Delta\psi+V\psi\leq 0,\quad\hbox{in}\ \Omega,\cr
           &{\partial\psi\over\partial\nu}\geq h\psi,\quad\quad\quad\hbox{on}\ \Sigma,\cr}\right.$$
and $v\in C^2(\Omega)\cap C^1(\bar\Omega)$ satisfies
$$\left\{\eqalign{&\Delta v+Vv\leq 0,\,\quad\hbox{in}\ \Omega,\cr
&{\partial v\over\partial\nu}\geq hv,\ \ \ \quad\quad\hbox{on}\ \Sigma,\cr
&v\geq 0,\ \,\qquad\qquad\hbox{on}\ \Gamma.\cr}\right.$$
Then $v\geq 0$ in $\bar\Omega$.
\medskip\rm
We now derive a Pohozaev-type identity for our problem; its proof is quite standard (see [\liuno]).
\medskip\noindent
{\bf Lemma 6.3.}\quad \it
Let $v$ be a $C^2$-solution of
$$\left\{\eqalign{& -\Delta v=0,\ \ \quad\quad\hbox{in}\ B^+_r,\cr
&{\partial v\over\partial\nu}=c(n)hv^q,\quad\hbox{on}\ \Gamma_1(B_r^+)=\partial B^+_r\cap\partial\erre^n_+,\cr}\right.\eqno\seitre$$
where $1\leq q\leq {n\over n-2}$ and $c(n)$ is constant depending on $n$. Then
$$\displaylines{c(n)\left({n-1\over q-1}-{n-2\over 2}\right)\int_{\Gamma_1(B_r^+)}hv^{q+1}\,d\sigma+{c(n)\over q+1}\int_{\Gamma_1(B_r^+)}\sum_{i=1}^{n-1}v^{q+1}{\partial h\over\partial x_i}x_i\, d\sigma\cr
-{c(n)r\over q+1}\int_{\partial\Gamma_1(B_r^+)}v^{q+1}h\,d\sigma'=\int_{\Gamma_2(B_r^+)}B(x,r,v,\nabla v)\, d\sigma\cr}$$
where $\Gamma_2(B_r^+)=\partial B^+_r\cap\erre^n_+$ and
$$B(x,r,v,\nabla v)={n-2\over 2}\,{\partial v\over\partial\nu}v+\mez r\left({\partial v\over\partial\nu}\right)^2-\mez r|\nabla_{\rm tan}v|^2$$
where $\nabla_{\rm tan}v$ denotes the component of the gradient $\nabla v$ which is tangent to $\Gamma_2(B_r^+)$.
\medskip\noindent
\rm An easy consequence of the previous Lemma is the following
\medskip\noindent
{\bf Corollary 6.4.}\quad \it Let $v(x)=a|x|^{2-n}+b+O(|x|)$ for $x$ close to $0$, with $a>0$ and $b>0$. There holds
$$\lim_{r\to 0^+}\int_{\Gamma_2(B_r^+)}B(x,v,\nabla v)<0.$$
\medskip\noindent\rm
In the proof of Lemma 4.1 we used the following result
\medskip\noindent
{\bf Lemma 6.5.}\quad\it Let $(M,g)$ be a smooth compact Riemannian manifold of positive type (namely $\lambda_1(B)>0$). Let $\eps_0>0$, $1\leq q\leq{n\over n-2}-\eps_0$. Suppose that $u$ satisfies
$$\left\{\eqalign{&-L_g u=0,\quad u>0,\quad\hbox{in}\ M,\cr
&{\partial u\over\partial \nu}=\mu u^q,\ \qquad\qquad\quad\hbox{on}\ \partial M,\cr
&\int_{\partial M}u^{q+1}=1.\cr}\right.\eqno\seinove$$
Then 
$$0<\mu=\int_M|\nabla_gu|^2+{n-2\over 4(n-1)}R_gu^2\leq C(M,g,\eps_0).$$
\medskip\noindent
{\bf Proof.}\quad \rm For $1+\eps_0\leq q\leq{n\over n-2}-\eps_0$, it follows from Theorem 3.1 that $C^{-1}\leq\mu^{1\over q-1}u\leq C$, which, together with $\int_{\partial M}u^{q+1}=1$, gives the claimed estimate. So we have to only establish the estimate for $1\leq q\leq 1+\eps_0$. We give a proof for $1\leq q\leq{n\over n-2}-\eps_0$. We can choose $f_1$ such that $E(f_1)=1$ and recall that $f_1$ satisfies
$$\left\{\eqalign{&-L_g f_1=0, \quad f_1>0,\,\quad\hbox{in}\ M,\cr
&{\partial f_1\over\partial \nu}=\lambda_1(B)f_1,\!\ \ \ \quad\quad\quad\hbox{on}\ \partial M.\cr}\right.$$
Multiply equation \seinove\ by $f_1$ and integrate by parts to obtain 
$$\mu\int_{\partial M}u^qf_1=\lambda_1(B)\int_{\partial M}f_1u\eqno\seidieci$$
which implies $\mu>0$. Note that, for $q=1$, $\mu=\lambda_1(B)$. In the following we assume $1<q<{n\over n-2}-\eps_0$. Since $1/c\leq f_1\leq c$ for some positive $c$, from \seidieci\  and the H\"older inequality, we deduce that
$$\mu\|u\|^{q-1}_{L^q(\partial M)}\leq c.\eqno\seiundici$$
From well-known interpolation inequalities, we deduce
$$\|u\|_{L^{q+1}(\partial M)}\leq\|u\|^\teta_{L^q(\partial M)}\|u\|^{1-\teta}_{L^{2(n-1)/(n-2)}(\partial M)}$$
where
$$\teta={q\over q+1}\cdot{n-nq-2q\over 2(n-1)-nq+2q}.$$
It is easy to check that $0<\teta<1$, $\teta^{-1}\leq c$, and $(1-\teta)^{-1}\leq c$. \par
Testing \seinove\ by $u$, we easily find that
$$\mu=\int_M\bigg(|\nabla_gu|^2+{n-2\over 4(n-1)}R_g u^2\bigg).$$
Therefore, from the Sobolev embedding Theorems, we deduce
$$1=\|u\|_{L^{q+1}(\partial M)}\leq c\|u\|^\teta_{L^q(\partial M)}\mu^{1-\teta\over 2}=c\left(\mu\|u\|_{L^q(\partial M)}^{2\teta\over 1-\teta}\right)^{1-\teta\over 2}.\eqno\seidodici$$
Combining \seiundici\ and \seidodici, we have that 
$$\mu^{1-{(1-\teta)(q-1)\over 2\teta}}\leq c.\eqno\seitredici$$
For $1\leq q\leq{n\over n-2}-\eps_0$, we have that
$$1-{(1-\teta)(q-1)\over 2\teta}\geq \delta(\eps_0)>0.\eqno\seiquattordici$$
The thesis follows from \seitredici\ and \seiquattordici.\hfill\fine
\medskip\noindent
The analogue for the negative case is
\medskip\noindent
{\bf Lemma 6.6.}\quad \it Let $(M,g)$ be a smooth compact Riemannian manifold with $\lambda_1(B)<0$ and $h_g\equiv 0$. Let $\eps_0>0$ and $1\leq q<\infty$. Suppose that $u$ satisfies \seinove. Then
$$0<-\mu=-\int_M\left(|\nabla_gu|^2+{n-2\over 4(n-1)}R_gu^2\right)\leq -{n-2\over 4(n-1)}\int_MR_gu^2\leq C(M,g).$$    

\goodbreak

\vskip2truecm
\centerline {\bf References}\bigskip\rm
\item {[\alm]} \maiuscolo A. Ambrosetti, Y. Y. Li, and A. Malchiodi: \sl Yamabe and scalar curvature problems under boundary conditions, \rm to appear, and preliminary notes on CRAS, s\'erie I \bf 330\rm(2000), 1013-1018.
 \smallskip
\item {[\changyang]} \maiuscolo S.A. Chang and P.C. Yang: \sl A perturbation result in prescribing scalar curvature on $S^n$, \rm Duke Mathematical Journal, \bf 64\rm (1991), 27-69.
\smallskip 
\item {[\cherrier]} \maiuscolo P. Cherrier: \sl Probl\`emes de Neumann non lin\'eaires sur les vari\'et\'es Riemanniennes, \rm J. Funct. Anal., \bf 57\rm(1984), 154-206.
 \smallskip 
\item {[\dma]} \maiuscolo Z. Djadli, A. Malchiodi, and M. Ould Ahmedou: \sl Prescribing scalar and mean curvature on the three dimensional half sphere, \rm preprint(2001).
\smallskip
 \item {[\escobaruno]} \maiuscolo J. F. Escobar: \sl Conformal deformation of a Riemannian metric to a scalar flat metric with constant mean curvature on the boundary, \rm Ann. of Math., (2) \bf 136\rm(1992), 1-50.
\smallskip 
\item {[\escobartre]} \maiuscolo J. F. Escobar: \sl The Yamabe problem on manifolds with boundary, \rm J. Diff. Geom., \bf 35\rm(1992), 21-84.
 \smallskip 
\item {[\escobardue]} \maiuscolo J. F. Escobar: \sl Conformal metrics with prescribed mean curvature on the boundary, \rm Calc. Var. and PDE's, \bf 4\rm(1996), 559-592. 
 \smallskip  
  \item {[\gilbarg]} \maiuscolo D. Gilbarg and N. S. Trudinger :
\rm ``Elliptic partial differential equations of second order'', second edition, Berlin, Springer-Verlag, 1983.
\smallskip
\item {[\hamza]} \maiuscolo H. Hamza: \sl Sur les transformations conformes des vari\'et\'es Riemanniennes a bord, \rm J. Funct. Anal., \bf 92\rm(1990), 403-447.
\smallskip 
 \item {[\hanliuno]} \maiuscolo Z. C. Han and Y. Y. Li:
\it The Yamabe problem on manifolds with boundary: exi\-sten\-ce and compactness results, \rm Duke Math. J., \bf  99\rm(1999), 489-542.
 \smallskip 
 \item {[\hanlidue]} \maiuscolo Z. C. Han and Y. Y. Li:
\sl The existence of conformal metrics with constant scalar curvature and constant boundary mean curvature, \rm Comm. Anal. Geom., \bf 8\rm(2000), 809-869.
 \smallskip 
 \item {[\hu]} \maiuscolo B. Hu:
\sl Nonexistence of a positive solution of the Laplace equation with nonlinear boundary condition, \rm Differential and integral equations, \bf 7\rm(1994), 301-313.
\smallskip 
 \item {[\kellogg]} \maiuscolo O. D. Kellog:
\rm ``Foundations of Potential Theory, \rm Dover Edition, 1954.
 \smallskip 
  \item {[\liuno]} \maiuscolo Y. Y. Li:
\sl Prescribing scalar curvature on $S^n$ and related topics, Part I, \rm J.
Differential Equations, \bf  120\rm(1995), 319-410.
\smallskip 
  \item {[\lidue]} \maiuscolo Y. Y. Li:
\sl The Nirenberg problem in a domain with boundary, \rm Topological methods in Nonlinear Analysis, \bf  6\rm(1995), 309-329.
% \smallskip 
% \item {[\lidue]} \maiuscolo Y. Y. Li: \sl Prescribing scalar curvature on $S^n$ and
%related topics, Part II: Existence and compactness, \rm Comm. Pure Appl.
%Math., \bf  49\rm(1996), 437-477.
\smallskip
\item {[\liliu]} \maiuscolo P. L. Li and J. Q. Liu: \sl Nirenberg's problem on the $2\!$-dimensional hemisphere, \rm Int. J. Math., \bf 4\rm(1993), 927-939.
 \smallskip 
\item {[\lizhuuno]} \maiuscolo Y. Y. Li and M. J. Zhu: \sl Uniqueness Theorems through the method of moving spheres, \rm Duke Math. J., \bf 80\rm(1995), 383-417.
 \smallskip 
\item {[\lizhudue]} \maiuscolo Y. Y. Li and M. J. Zhu: \sl Yamabe type equations on three dimensional Riemannian manifolds, \rm Communications in contemporary Math., \bf 1\rm(1999), 1-50.
 \smallskip 
\item {[\louzhu]} \maiuscolo Y. Lou and M. J. Zhu: \sl Classification of nonnegative solutions to some elliptic problems, \rm Diff. Integral Eq., \bf 12\rm(1999), 601-612.
 \smallskip
\item {[\ma]} \maiuscolo L. Ma: \sl Conformal metrics with prescribed mean curvature on the boundary of the unit ball, \rm to appear.
 \smallskip  
  \item {[\nirenberg]} \maiuscolo L. Nirenberg:
\rm ``Topics in nonlinear functional analysis'', \rm lectures notes, 1973-1974, Courant Institute of Mathematical Sciences, NewYork University, New York, 1974.
 \smallskip
\item {[\ahmedou]} \maiuscolo M. Ould Ahmedou: \sl A Riemann mapping Theorem in higher dimensions, Part I: the conformally flat case with umbilic boundary, \rm preprint(2001).
 \smallskip
\item {[\schoen]} \maiuscolo R. Schoen: \rm ``Variational theory for the total scalar curvature functional for Riemannian metrics and related topics'' in \sl Topics in Calculus of Variations (Montecatini Terme, 1987), \rm Lectures Notes in Math., \bf 1365\rm, Springer-Verlag, Berlin, 1989, 120-154.
 \smallskip 
 \item {[\schoenyau]} \maiuscolo R. Schoen and S. T. Yau:
\sl Conformally flat manifolds, Kleinian groups and scalar curvature, \rm Invent. Math., \bf  92\rm(1988), 47-71.
 \smallskip 
\smallskip
 \item {[\schoenzhang]} \maiuscolo R. Schoen and D. Zhang:
\sl Prescribed scalar curvature on the $n\!$-sphere, \rm Calc. Var. and PDE's, \bf  4\rm(1996), 1-25.
 \smallskip

\end

In order to apply Theorem 6.2 we have to check that
$$\leqalignno{&\Delta\varphi_i\leq0,\qquad\qquad\quad\qquad\qquad\hbox{in}\ D_i,&{\rm(i)}\cr
            &{\partial\varphi_i\over\partial x^n}+(n-2)f_i^{\tau_i}\vi^{\qi-1}\leq	0,\quad\hbox{on}\ \Gamma_1(D_i),\hskip6truecm&{\rm (ii)}\cr
            &\varphi_i\geq0,\ \ \,\qquad\qquad\qquad\qquad\quad\hbox{on}\  \Gamma_2(D_i).&{\rm(iii)}\cr}$$

%PROOF OF LEMMA 2.8
\medskip\noindent
{\bf Proof.}\quad\rm From Theorem 6.3 we have that for $0<r<1$
$$\displaylines{(n-2)\left({n-1\over\qi-1}-{n-2\over 2}\right)\int_{\Gamma_1(B_r^+)}f_i^{\ti}\vi^{\qi+1}\,d\sigma+{n-2\over \qi-1}+\int_{\Gamma_1(B_r^+)}\sum_{j=1}^{n-1}\vi^{\qi+1}{\partial(f_i^{\ti})\over\partial x_j}x_j\,d\sigma\cr
 -{(n-2)r\over\qi+1}\int_{\partial\Gamma_1(B_r^+)}\vi^{\qi+1}f_i^{\ti}\,d\sigma=\int_{\Gamma_2(B_r^+)}
B(x,r,\vi,\nabla\vi)\,d\sigma\cr}$$
where $B$ is given by
$$B(x,r,\vi,\nabla\vi)={n-2\over2}{\partial \vi\over\partial\nu}\vi+\mez r\left({\partial \vi\over\partial\nu}\right)^2-\mez r|\nabla_{\rm tan}\vi|^2\eqno\defdibi$$
where $\nabla_{\rm tan}\vi$ is the tangent component of $\nabla\vi$. From Lemma 2.6, Lemma 2.7, and some standard elliptic estimates we derive the claimed estimate.\hfill\fine

%PROOF OF POHOZAEV
{\bf Proof.}\quad\rm Multiply \seitre\ by $v$ and integrate by parts over $B^+_r$
$$\eqalign{0=\int_{B_r^+}(-\Delta v)v\, dx&=\int_{B_r^+}|\nabla v|^2\, dx-\int_{\partial B_r^+}{\partial v\over\partial\nu}v\, d\sigma\cr
                   &=\int_{B_r^+}|\nabla v|^2\, dx-\int_{\Gamma_1(B_r^+)}{\partial v\over\partial\nu}v\, d\sigma-\int_{\Gamma_2(B_r^+)}{\partial v\over\partial\nu}v\, d\sigma.\cr}$$
Hence, in view of the boundary condition in \seitre, we obtain 
$$\int_{B_r^+}|\nabla v|^2\, dx=\int_{\Gamma_1(B_r^+)}c(n)hv^{q+1}\, d\sigma+\int_{\Gamma_2(B_r^+)}{\partial v\over\partial\nu}v\, d\sigma.\eqno\seiquattro$$
Multiply \seitre\ by $x\cdot\nabla v$ and integrate by parts over $B_r^+$
$$\displaylines{0=\int_{B_r^+}(-\Delta v)\bigg(\sum_{i=1}^{n}{\partial v\over\partial x_i}x_i\bigg)\,dx=\int_{B_r^+}\nabla v\cdot\nabla\bigg(\sum_{i=1}^{n}{\partial v\over\partial x_i}x_i\bigg)\,dx\cr
\quad-\int_{\Gamma_1(B_r^+)}{\partial v\over\partial\nu}\bigg(\sum_{i=1}^{n}{\partial v\over\partial x_i}x_i\bigg)\,d\sigma-\int_{\Gamma_2(B_r^+)}{\partial v\over\partial\nu}\bigg(\sum_{i=1}^{n}{\partial v\over\partial x_i}x_i\bigg)\,d\sigma.\cr}$$
Hence, in view of the boundary condition in \seitre, we have that
$$I+I\!\!I+I\!\!I\!\!I=0$$
where
$$\eqalign{I&=\int_{B_r^+}\nabla v\cdot\nabla\bigg(\sum_{i=1}^{n}{\partial v\over\partial x_i}x_i\bigg)\,dx,\cr
I\!\!I&=-\sum_{i=1}^n\int_{\Gamma_1(B_r^+)}c(n)hv^q{\partial v\over\partial x_i}x_i\, d\sigma,\cr
I\!\!I\!\!I&=-\int_{\Gamma_2(B_r^+)}{\partial v\over\partial\nu}\bigg(\sum_{i=1}^{n}{\partial v\over\partial x_i}x_i\bigg)\,d\sigma.\cr}$$
After an integration by parts in $I\!\!I$, we have that
$$\eqalign{I\!\!I&=-{1\over q+1}\sum_{i=1}^{n-1}c(n)\int_{\Gamma_1(B_r^+)}hx_i{\partial\over\partial x_i}(v^{q+1})\,d\sigma\cr
&={1\over q+1}c(n)\sum_{i=1}^{n-1}\int_{\Gamma_1(B_r^+)}v^{q+1}{\partial\over\partial x_i}(hx_i)\,d\sigma-{1\over q+1}c(n)\sum_{i=1}^{n-1}\int_{\partial\Gamma_1(B_r^+)}v^{q+1}hx_i\nu_i\,d\sigma'.\cr}$$
Hence, using the fact that on $\partial\Gamma_1(B_r^+)$ $x\cdot\nu=|x|=r$, we derive
$$\eqalignno{I\!\!I=\,&c(n){n-1\over q+1}\int_{\Gamma_1(B_r^+)}hv^{q+1}\,d\sigma+{c(n)\over q+1}\sum_{i=1}^{n-1}\int_{\Gamma_1(B_r^+)}v^{q+1}{\partial h\over\partial x_i}x_i\,d\sigma&\cr
&-{c(n)r\over q+1}\int_{\partial \Gamma_1(B_r^+)}v^{q+1}h\,d\sigma'.&\seicinque\cr}$$
As far as $I\!\!I\!\!I$ is concerned, since on $\Gamma_2(B_r^+)$ $x=r\nu$, we have that
$$I\!\!I\!\!I=-r\int_{\Gamma_2(B_r^+)}{\partial v\over\partial\nu}\nabla v\cdot\nu\,d\sigma=-r\int_{\Gamma_2(B_r^+)}\left({\partial v\over\partial\nu}\right)^2\,d\sigma.\eqno\seisei$$
Let us now consider $I$:
$$\eqalign{I&=\int_{B_r^+}\sum_{j=1}^n{\partial v\over\partial x_j}{\partial\over \partial x_j}\bigg(\sum_{i=1}^n{\partial v\over\partial x_i}x_i\bigg)\,dx\cr
&=\int_{B_r^+}\sum_{i,j=1}^n{\partial v\over\partial x_j}{\partial v\over \partial x_i}\delta_{i,j}\,dx+\int_{B_r^+}\sum_{i,j=1}^n{\partial v\over\partial x_j}{\partial\over \partial x_j}\bigg({\partial v\over\partial x_i}\bigg)x_i\,dx\cr
&=\int_{B_r^+}|\nabla v|^2\,dx+\mez\sum_{i,j=1}^n\int_{B_r^+}{\partial \over\partial x_i}\bigg(\bigg|{\partial v\over\partial x_j}\bigg|^2\bigg)x_i\,dx\cr
&=\int_{B_r^+}|\nabla v|^2\,dx-{n\over 2}\sum_{j=1}^n\int_{B_r^+}\bigg|{\partial v\over\partial x_j}\bigg|^2\,dx+\mez\sum_{i,j=1}^n\int_{\partial B_r^+}\bigg|{\partial v\over\partial x_j}\bigg|^2x_i\nu_i\,d\sigma\cr
&={2-n\over 2}\int_{B_r^+}|\nabla v|^2\,dx+\mez\sum_{j=1}^n\int_{\partial B_r^+}\bigg|{\partial v\over\partial x_j}\bigg|^2x\cdot\nu\,d\sigma.\cr}$$
Hence, after noting that on $\Gamma_1(B_r^+)$ $x\cdot\nu=0$, we obtain
$$I={2-n\over 2}\int_{B^+_r}|\nabla v|^2\,dx+\mez\sum_{j=1}^n\int_{\Gamma_2(B_r^+)}\bigg|{\partial v\over\partial x_j}\bigg|^2 x\cdot \nu\,d\sigma.\eqno\seisette$$
From \seiquattro\ and \seisette\ we deduce that
$$\eqalignno{I=&{2-n\over 2}\int_{\Gamma_1(B_r^+)}c(n)hv^{q+1}\,d\sigma&\cr
&+{2-n\over2}\int_{\Gamma_2(B_r^+)}{\partial v\over\partial\nu}v\,d\sigma+\mez\sum_{j=1}^n\int_{\Gamma_2(B_r^+)}\bigg|{\partial v\over\partial x_j}\bigg|^2 x\cdot \nu\,d\sigma.&\seiotto\cr}$$
Putting together \seicinque, \seisei, and \seiotto, from $I+I\!\!I+I\!\!I\!\!I=0$, we deduce that
$$\displaylines{c(n)\left({n-1\over q+1}-{n-2\over2}\right)\int_{\Gamma_1(B_r^+)}hv^{q+1}+{c(n)\over q+1}\int_{\Gamma_1(B_r^+)}\sum_{i=1}^{n-1}v^{q+1}{\partial h\over\partial x_i}x_i\,d\sigma\cr
-{c(n)r\over q+1}\int_{\partial\Gamma_1(B_r^+)}v^{q+1}h\,d\sigma'=\int_{\Gamma_2(B_r^+)}B(r,x,v,\nabla v)\,d\sigma,\cr}$$
where
$$B(r,x,v,\nabla v)={n-2\over2}\,{\partial v\over\partial\nu}v-\mez\sum_{j=1}^n\left|{\partial v\over\partial x_j}\right|^2 x\cdot\nu+r\left({\partial v\over\partial\nu}\right)^2.$$
On $\Gamma_2(B_r^+)$, $x\cdot \nu=r$ and so
$$|\nabla v|^2 x\cdot\nu=\left({\partial v\over\partial\nu}\right)^2x\cdot\nu+|\nabla_{\rm tan}v|^2x\cdot \nu=r\left({\partial v\over\partial\nu}\right)^2+r|\nabla_{\rm tan}v|^2.$$
Then, on $\Gamma_2(B_r^+)$, we can write $B$ in the form
$$B(r,x,v,\nabla v)={n-2\over2}\,{\partial v\over\partial\nu}v+\mez r\left({\partial v\over\partial \nu}\right)^2-\mez r|\nabla_{\rm tan}v|^2.$$
The Lemma is proved.\hfill\fine